\documentclass[12pt]{article}
\usepackage{amsmath}
\usepackage{amssymb}
\usepackage{amsthm}
\usepackage{latexsym}
\usepackage{color}
\usepackage{graphicx}
\usepackage{color}
\usepackage{caption}
\usepackage{subcaption}
\usepackage[utf8]{inputenc} 
\usepackage{float}

\DeclareSymbolFont{calletters}{OMS}{cmsy}{m}{n}
\DeclareSymbolFontAlphabet{\mathcal}{calletters}

\def\be{\begin{align}}
\def\ee{\end{align}}
\def\bef{\begin{flalign}}
\def\eef{\end{flalign}}

\newtheorem{Theorem}{Theorem}[section]

\newtheorem{Remark}[Theorem]{Remark}

\makeatletter \@addtoreset{equation}{section}

\newcommand{\un}{1\hspace{-1mm}{\rm I}}   

\usepackage{color}

\def \E{\mathbb{E}}

\def \R{\mathbb{R}}
\def \S{\mathbb{S}}
\def \M{\mathbb{M}}
\def \N{\mathbb{N}}

\def\Kc{{\cal K}}
\def\Lc{{\cal L}}

\def\Wc{{\cal W}}

\def\Kct{\widetilde{\Kc}}

\def \I{{\bf I}}

\def \0{\mathbf{0}}

\def \Kcb{\overline{\Kc}}
\def \psih{\widehat{\psi}}

\def\x{\times}

\def\1{{\bf 1}}

\def \Fb{\overline F}

\def \Kch{\widehat \Kc}

\usepackage{mathtools}

\title{Variations on branching methods for non linear PDEs}

\author{ Xavier Warin \thanks{EDF R\&D \& FiME, Laboratoire de Finance des March\'es de l'Energie, xavier.warin@edf.fr}}

\date{\today}

\begin{document}

\maketitle

\abstract{ The branching methods developed in \cite{LOTTW}, \cite{LTTW} are effective methods to solve some  semi linear PDEs and are shown
numerically to be  able to solve some full non linear PDEs.
These  methods are however restricted to some small coefficients in the PDE and small maturities.
This article shows numerically that these methods can be adapted to solve the problems with longer maturities in the semi-linear case by using a new derivation scheme and some nested method.
As for the case of full non linear PDEs, we introduce new schemes  and we show numerically that they provide an effective alternative to the schemes previously developed.
}

\section{Introduction}
The resolution of low dimensional non linear PDEs  is often achieved by some deterministic methods such as finite difference schemes, finite elements and finite volume.
Due the curse of dimensionality, these methods cannot be used in dimension greater than three : both the computer time and the memory required are too large even for supercomputers.
In the recent years the probabilistic community has developed some representation of semi linear PDE:
\begin{flalign}
\label{eq:semiLin}
 -\partial_tu-\Lc u =f(u,Du), & \quad \quad \quad \quad u_T=g,& t<T,~x\in\R^d,
\end{flalign}
by means of backward stochastic differential equations (BSDE), as introduced by  \cite{PardouxPeng}.
Numerical Monte Carlo algorithms have been developed to solve efficiently these BSDE by  \cite{BouchardTouzi},  \cite{zhang}.
The representation  of  the following full non linear PDE: 
\begin{flalign}
 \label{eq:nonLin}
-\partial_tu-\Lc u =f(u,Du,D^2u), & \quad \quad \quad \quad u_T=g,& t<T,~x\in\R^d,
\end{flalign}
has been given by the mean of second order backward stochastic differential  equation (SOBSDE) by  \cite{cheridito}.
A numerical algorithm  developed by \cite{FTW} has been derived  to solve these  full non linear PDE by the mean of SOBSEs.\\
The BSDE and SOBSDE schemes developed rely on the approximation of conditional expectation and the most effective implementation is based on regression methods as developed in  \cite{LGW1}, \cite{LGW2}.
These regression methods develop an approximation of   conditional expectations based on an expansion  on basis functions. The size of this expansion has to grow exponentially with the dimension of the problem so we have to face again the curse of dimensionality. 
Notice that the BSDE methodology could be used in dimension 4 or 5 as regressions has been successfully used in dimension 6 in   \cite{BW} using some local regression function.\\
Recently a new representation of semi linear equations \eqref{eq:semiLin} for a polynomial function $f$ of $u$ and $Du$ has been given by \cite{LOTTW} : this representation
uses the automatic differentiation approximation  as used in \cite{FLLLT}, \cite{BET}, \cite{HTT2}, and \cite{DOW}. 
The authors have shown that the representation gives a finite variance estimator only for small maturities or small non linearities and numerical  examples until dimension 10 are given. Besides, they have shown that the given scheme using Malliavin weights cannot be used to solve the full non linear equation \eqref{eq:nonLin}.
\\
\cite{LTTW} have introduced a re-normalization technique improving numerically the convergence of the scheme diminishing the variance observed for the semi linear case. Besides, the authors haved introduced a scheme to solve the full non linear equation \eqref{eq:nonLin}. Without proof of convergence they numerically have shown that the developed scheme is effective.\\
The aim of the paper is to provide some numerical variation on the algorithm developed in  \cite{LOTTW,LTTW}.
In a first part we will show, with simple ideas, that it is possible to deal with longer maturities than the ones possible with the initial algorithm.\\
In a second part we give some alternative schemes to the one proposed in \cite{LTTW} and, testing  them on some numerical examples, we show that they are superior than the scheme previously developed.\\
In the numerical results presented in the article, all errors are estimated  as  the $\log$ of the  standard deviation observed divided  by the square root of the number of particles used and these errors are plotted as a function of the $\log$ of the number of particles used. As our methods are pure Monte Carlo methods we expect to have lines with slope $-\frac{1}{2}$ when the numerical variance is bounded.

\section{The Semi Linear case}
Let $\sigma_0 \in \S^d$ be some constant non-degenerate matrix,  $\mu \in \R^d$ be some constant vector  $f: [0,T] \x \R^d \x \R \x \R^d \to \R$,
and $g: \R^d \to \R$ bounded Lipschitz functions, 
we consider the semi linear parabolic PDE:
\begin{flalign} \label{eq:PDE}
\partial_t u + \frac{1}{2} \sigma_0 \sigma_0^{\top} : D^2 u + \mu.Du + f(\cdot, u ,Du) = 0,
~~~\mbox{on}~~[0,T) \x \R^d,
\end{flalign}
with terminal condition $u(T, \cdot) = g(\cdot)$ where 	$A:B :=  \mbox{Trace}(AB^{\top})$ for two matrices $A$, $B \in \M^d$.\\
When $f$ is a polynomial in $(u, Du)$ in the form
\begin{flalign*}
f(t,x,y,z,\gamma)
~=\!
\sum_{\ell = (\ell_0, \ell_1, \cdot, \ell_m) \in L}
\! c_{\ell}(t,x) y^{\ell_0} \prod_{i=1}^m (b_i \cdot z)^{\ell_i} ,
\end{flalign*}
for some $m \ge 1$, $L\subset \N^{1 + m}$, 
where $(b_i)_{i=1, m}$ is a sequence of $\R^d-$valued bounded continuous functions defined on $[0,T] \x \R^d$, 
and $(c_{\ell})_{\ell \in L}$ is a sequence of bounded continuous functions defined on $[0,T] \x \R^d$.
 \cite{LOTTW} obtained a probabilistic representation to the above PDE by branching diffusion processes under some technical conditions.
In the sequel, we simplify the setting by taking $f$  as a constant (in $u$, $Du$) plus a monomial  in $u$, $(b_i.Du)$ , $i=1,m$ : 
\begin{flalign} \label{eq:generator}
f(t,x,y,z)
~=\!
 h(t,x)+ c(t,x) y^{\ell_0} \prod_{i=1}^m (b_i \cdot z)^{\ell_i} ,
\end{flalign}
for some $m \ge 1$,
where $(b_i)_{i=1, m}$  is a sequence of $\R^d-$valued bounded continuous function defined on $[0,T] \x \R^d$, $(\ell_i)_{i=0,m} \in \N^{m+1}$ supposing that $\displaystyle{\sum_{i=0,m}}\ell_i >0$,
and $c$ is a  bounded continuous function defined on $[0,T] \x \R^d$. We note $L= \sum_{i=0}^m \ell_i$.
\begin{Remark}
The case with $f$ a general polynomial only complexifies the notation : it can be simply treated as in \cite{LOTTW} by introducing some probability mass function $(p_{\ell})_{\ell \in L}$ (i.e. $p_{\ell} \ge 0$ and $\sum_{\ell \in L} p_{\ell} = 1$) that are used to select with monomial to consider during the branching procedure. Another approach can be used : instead of sampling the monomial to use, it is possible to consider successively  all terms of the $f$ but this doesn't give a representation as nice as the one in  \cite{LOTTW}.
\end{Remark}
\subsection{Variation on the original scheme of \cite{LOTTW}}
\label{origScheme}
In this section we present the original scheme of \cite{LOTTW} and explain how to diminish the variance increase the maturities of the problem.
\subsubsection{The branching process}
Let us first introduce a branching process with arrival time of distribution density function $\rho$.  At the arrival time, the particle branches into $|\ell |$ offsprings.
We introduce a sequence  of   i.i.d. positive random variables $(\tau^{k})_{k = (k_1, \cdots, k_{n-1}, k_n) \in \N^n, n>1}$ with all the values $k_i \in [1,L]$, for  $i >0$.\\
We construct an age-dependent branching process using the following procedure :
\begin{enumerate}
\item  We start from a particle marked by $0$, indexed by $(1)$, of generation $1$,
		whose arrival time is given by $T_{(1)} := \tau^{(1)} \wedge T$.
\item  Let $k = (k_1, \cdots, k_{n-1}, k_n) \in \N^n$ be a particle of generation $n$, with arrival time $T_k$ that branches into $L$ offspring particles noted $(k_1, \cdots, k_{n-1}, k_n,i)$ for $i=1,...,L$. We  define the set of its offspring particles by
	$$S(k) := \{(k_1, \cdots, k_n, 1), \cdots, (k_1, \cdots, k_n, L) \},$$
  We first mark  the $\ell_0$ particles by 0,  the $\ell_1$ next by 1 , and so on,  so that each particle has a mark $i$ for $i = 0, \cdots, m$.
\item For a particle $k = (k_1, \cdots, k_n, k_{n+1})$ of generation $n+1$, 
  we denote by $k- := (k_1, \cdots, k_n)$ the ``parent'' particle of $k$,
  and the arrival time of $k$ is given by $T_k := \big(T_{k-} + \tau^{k} \big) \wedge T$. Let us denote $ \Delta T_k = T_k -T_{k-}$.
\item In particular, for a particle $k = (k_1, \cdots, k_n)$ of generation $n$,
  and $T_{k-}$ is its birth time and also the arrival time of $k-$.
  Moreover, for the initial particle $k = (1)$, one has $k- = \emptyset$, and $T_{\emptyset} = 0$.
\end{enumerate}
We denote further
$$
\theta_k := \mbox{mark of}~k,
~~~
\Kc^n_t := \begin{cases}
  \big\{ k ~\mbox{of generation}~n~\mbox{s.t.}~T_{k-} \le t < T_k \big\}, &\mbox{when}~~t \in [0,T),\\
  \{k ~\mbox{of generation}~ n~\mbox{s.t.}~ T_k = T\}, &\mbox{when}~~ t = T,
\end{cases}
$$
and also
$$
\Kcb^n_t := \cup_{s \le t} \Kc^n_s,
~~~~
\Kc_t := \cup_{n \ge 1} \Kc^n_t
~~~\mbox{and}~~~~
\Kcb_t := \cup_{n \ge 1} \Kcb^n_t.
$$
Clearly, $\Kc_t$ (resp. $\Kc^n_t$) denotes the set of all living particles (resp. of generation $n$) in the system at time $t$,
and $\Kcb_t$ (resp. $\Kcb^n_t$) denotes the set of all particles (resp. of generation $n$) being alive at or before time $t$.\\

We next equip each particle with a Brownian motion in order to define a branching Brownian motion.
Let $(\hat W^{k})_{k = (k_1, \cdots, k_{n-1}, k_n) \in \N^n, n>1}$ be a sequence of independent $d$-dimensional Brownian motion, which is also independent of $(\tau^{k})_{k = (k_1, \cdots, k_{n-1}, k_n) \in \N^n, n>1}$.
Define $W^{(1)}_t =  \hat W^{(1)}_t$ for all $t \in \big[0, T_{(1)} \big]$ and then for each $k = (k_1, \cdots, k_n) \in \Kcb_T \setminus \{(1)\}$,
define
\begin{flalign}\label{eq:def_Wk}
W^k_t ~:=~ W^{k-}_{T_{k-}} + \hat  W^k_{t - T_{k-}}, ~~\mbox{for all}~ t \in [T_{k-}, T_k].
\end{flalign}
Then $(W^k_{\cdot})_{k \in \Kcb_T}$ is a branching Brownian motion.
\subsubsection{The original algorithm}
\label{sec:origAlgo}
Let us denote  $\bar F(t):=\int_t^\infty\rho(s)ds$.
Denoting $X^k_t := x + \mu t + \sigma_0 W^k_t$ for all $k \in \Kcb_T$ and $t \in [T_{k-}, T_k]$ and
by $\E_{t,x}$ the expectation operator conditional on the starting data $X_t=x$ at time $t$, 
 we obtain from the Feynman-Kac formula the representation of the solution $u$ of equation \eqref{eq:PDE} as:
 \begin{flalign}
\label{eq:uVal}
 u(0,x)
 =
 \E_{0,x} \Big[\bar F(T)\frac{g(X_T)}{\bar F(T)}+\int_0^T \frac{f(u,Du)(t,X_t)}{\rho(t)}\rho(t)dt\Big]
 =
 \E_{0,x} \big[\phi\big(T_{(1)},X^{(1)}_{T_{(1)}}\big)\big],
 \end{flalign}
where $T_{(1)}:=\tau^{(1)}\wedge T$, and
  \begin{flalign}
  \phi(t,y)
  := \frac{\1_{\{t\ge T\}}}{\bar F(T)}g(y)
                    \!+\! \frac{\1_{\{t<T\}}}{\rho(t)}( h+c u^{\ell_0} \prod_{i=1}^m (b_i \cdot Du)^{\ell_i})(t,y).
 \label{phi}
 \end{flalign}
On the event $\{  \1_{\{T_{(1)}<T\}}\}$, using the independence of the $(\tau^k,W^k)$   we are left to calculate

\begin{flalign}
\label{eq:calU}
[c u^{\ell_0} & \displaystyle{ \prod_{i=1}^m (b_i \cdot Du)^{\ell_i}](T_{(1)},X_{T_{(1)}})  = c \prod_{j=1}^{\ell_0} \E_{T_{(1)},X_{T_{(1)}}} \big[ \phi\big(T_{(1,j)},X^{(1)}_{T_{(1,j)}}\big)\big]}  \nonumber  \\ & \displaystyle{ \prod_{i=1}^m ( b_i(T_{(1)},X_{T_{(1)}}).D \E_{T_{(1)},X_{T_{(1)}}}\big[\phi\big(T_{(1,p)},X^{(1,p)}_{T_{(1,p)}}\big)\big])^{\ell_i}}
\end{flalign}
Using differentiation with respect to the heat kernel, i.e. the marginal density of the Brownian motion we get :
\begin{flalign}
\label{eq:RecVal}
[c u^{\ell_0} \prod_{i=1}^m (b_i \cdot Du)^{\ell_i}](T_{(1)},X_{T_{(1)}})  = c \prod_{j=1}^{\ell_0} \E_{T_{(1)},X_{T_{(1)}}} \big[ \phi\big(T_{(1,j)},X^{(1,j)}_{T_{(1,j)}}\big)\big]  \nonumber \\ \quad \quad  \prod_{i=1}^m ( b_i(T_{(1)},X_{T_{(1)}}).\E_{T_{(1)},X_{T_{(1)}}}\big[ (\sigma_0^\top)^{-1}\frac{ \hat W^{(1,p)}_{\Delta T_{(1,p)}}}{\Delta T_{(1,p)}} \phi\big(T_{(1,p)},X^{(1,p)}_{T_{(1,p)}}\big)\big])^{\ell_i}
\end{flalign} 
Using equations \eqref{eq:uVal} and \eqref{eq:RecVal} recursively and the tower property , we get the following
representation 
\begin{flalign}
\label{eq:initRep}
u(0,x)
 = \E_{0,x}\Big[  \psih_{(1)}\Big]
\end{flalign}
where $\psih_{(1)}$ is given by the backward recursion :
let  $\psih_k := \frac{g(X^k_T) -g(X^k_{T_{k-}})\1_{\{\theta_k \neq 0\}}}{\Fb(\Delta T_k)}$ for every $k \in \Kc_T$, then let
\begin{flalign}\label{eq:backRep}
\psih_k
~:=~
\frac{1}{\rho(\Delta T_k)} \big( h(T_k,X^k_{T_k}) + c(T_k, X^k_{T_k})
\prod_{{\tilde k} \in S(k)} \!\!\! \psih_{{\tilde k}} \Wc_{{\tilde k}} \big),
~~~~\mbox{for}~k \in  \Kcb_T \setminus \Kc_T.
\end{flalign}

where 
\begin{flalign}
\label{eq:weigthSem}
		 \Wc_k
		 ~=~
		 \1_{\{\theta_k = 0\}}
		~+~
		\1_{\{\theta_k \neq 0\}}
		~\frac{ b_{\theta_k}(T_{k-}, X^k_{T_{k-}})
		\cdot (\sigma_0^{\top})^{-1} \hat W^k_{\Delta T_k}}
		{\Delta T_k}.
\end{flalign}
and we have used that  $\E_{0,x}\big[ g(X^k_{T_k-})  b_{\theta_k}(T_{k-}, X^k_{T_{k-}})
		\cdot \sigma_0^{\top})^{-1} \hat W^k_{\Delta T_k}\big]  =0$.
This backward representation is slightly different from the  elegant representation  introduced  in \cite{LOTTW}.
Clearly on our case the variance of the method used will be lower than with the representation   in \cite{LOTTW} for a similar computational cost.\\

In the case where the operator $f$ is linear  and a function of the gradient ($\ell_0=0$, $m=1$ and $\ell_1=1$) using the arguments in \cite{DOW}  it can be easily seen by conditioning with respect to the number of branching that equation \eqref{eq:initRep} is of finite variance if $\frac{1}{x \rho(x)^2} = O(x^\alpha)$ as $x \longrightarrow 0$ with $\alpha \ge 0$.\\
When $\tau$ follows for example a gamma law with parameters $\kappa$ and $\theta$, the finite variance is proved as soon as $\kappa \le 0.5$ for PDE coefficients and maturities small enough.\\
In the non linear case, \cite{LOTTW} have shown that the variance is in fact finite  for maturities small enough and small coefficients  as soon as $\kappa <0.5$ but numerical results show that $\kappa=0.5$ is optimal in term of efficiency: for  a given $\theta$ the numerical variance is nearly the same for the values of $\kappa$ between 0.4 and 0.5 but a higher $\kappa$ value limits the number of branching thus meaning a smaller computational cost.\\

\subsubsection{Variation on the original  scheme}
As indicated in the introduction, the method is restricted to small maturities or small non linearities.
Having a given non linearity we are interested in adapting the methodology in  order to be able to treat longer maturities. A simple idea consists in noting that the Monte Carlo  method is applied by sampling the conditional expectation $\E_{t,x}$ for $t>0$ appearing in equation \eqref{eq:RecVal} only once. Using nested Monte Carlo, so by sampling each term of equation \eqref{eq:RecVal} more that one time
one can expect a reduction in the variance observed. A nested method of order $n$ is defined as a method using $n$ sampling to estimate each function $u$ or $Du$ at each branching.
Of course the computational time will grow exponentially with the number of samples taken and for example 
trying to use a gamma law  with a non linearity of Burger's type  $u(b.Du)$  with $\kappa =0.5$  is very costly: due to the  high values of the density $\rho$ near $0$, trajectories can have many branching.\\
Some different strategies have been tested to  be able to use this technique :
\begin{itemize}
\item A first possibility consists in trying  to re-sample more at the beginning of the resolution and decreasing the number of samples as time goes by or as the number of branching increases.  The methodology works  slightly better than  a re-sampling with a constant number of particles but  has to be adapted to each  maturity and each case so it has been given up.
\item  Another observation is that the gamma law is only necessary  to treat the gradient term: so it is possible to use two laws: a first one, an exponential law,  will be used to estimate the $u$ function while an gamma law will be used for the $Du$ terms.  This second technique is the most effective and is used for the results obtained in the section.
\end{itemize}
For a given dimension $d$ , we take $\sigma_0 = \frac{1}{\sqrt{d}} \I_d$, $\mu =\0$,
\begin{flalign*}
 f(t,x,y,z)= d(t,x) +  y (b \cdot z),
\end{flalign*}
where $b := \frac{0.2}{d} (1+\frac{1}{d}, 1+\frac{2}{d}, \cdots, 2)$ and 
\begin{flalign*}
h(t,x) 
:=
\cos( x_1 +\cdots + x_d) 
\Big( \alpha + \frac{\sigma_0^2}{2}  + c \sin(x_1 +\cdots + x_d)  \frac{3d+1}{2d} e^{ \alpha (T-t)} \Big) 
e^{ \alpha (T-t)}.	
\end{flalign*}
With terminal condition $g(x) = \cos( x_1+ \cdots + x_d)$,
the explicit solution of semi linear PDE \eqref{eq:PDE} is given by
$$
u(t,x) =  \cos(x_1+ \cdots + x_d) e^{\alpha (T-t)}.
$$
Our goal is to estimate $u$ at $t=0$, $x=0.5 \1$.
This test case will be noted test A in the sequel. \\
We use the nested algorithm with two distributions for $\tau$:
\begin{itemize}
\item an exponential law with density $\rho(s)= \lambda e^{-\lambda s}$ with $\lambda=0.4$ to calculate the $u$ terms,
\item a gamma distribution $\rho(s) ~=~ \frac{1}{\Gamma(\kappa) \theta^{\kappa}} s^{\kappa -1} \exp(- s/\theta) \1_{\{s > 0\}}$ with \newline $\Gamma(\kappa) := \int_0^{\infty} s^{\kappa-1} e^{-s} ds$ and the parameters $\kappa=0.5$, $\frac{1}{\theta}=0.4$ to calculate the $Du$ terms.
\end{itemize}
 
\begin{figure}[H]
  \centering
  \includegraphics[width=7cm]{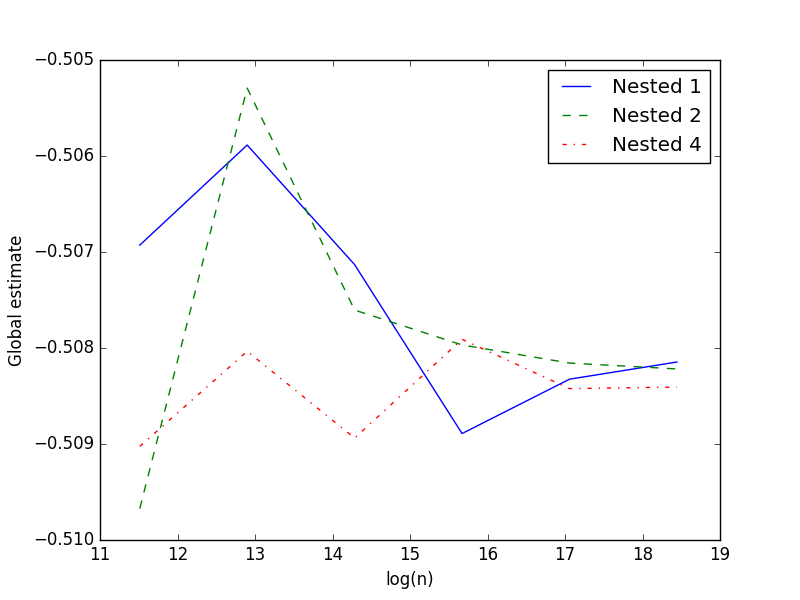}
  \includegraphics[width=7cm]{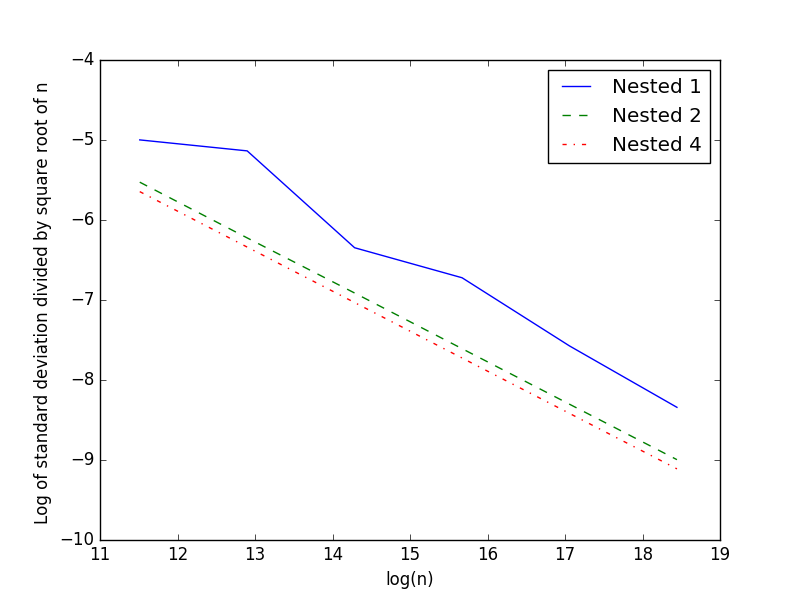}
  \caption{Estimation and error in $d=4$  on case test A. Maturity $T=1$. }
\label{semiFig1}
\end{figure}

\begin{figure}[H]
  \centering
  \includegraphics[width=7cm]{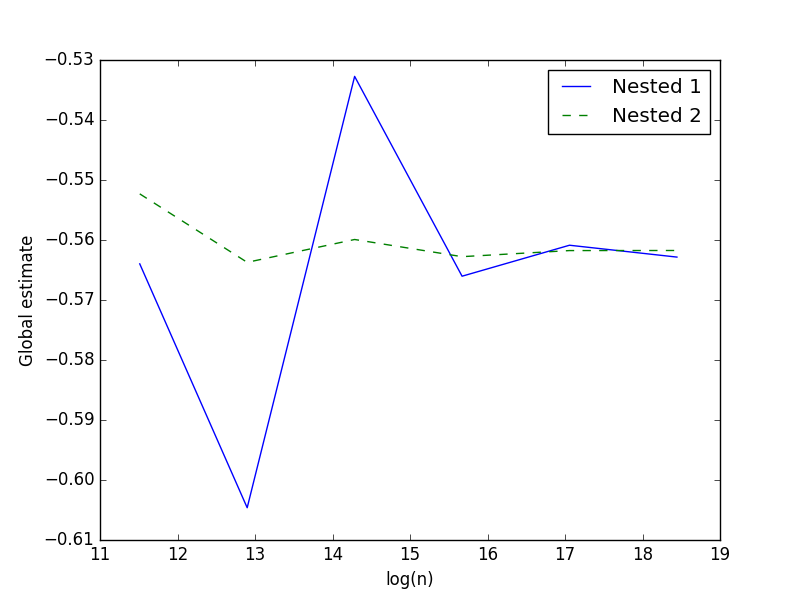}
  \includegraphics[width=7cm]{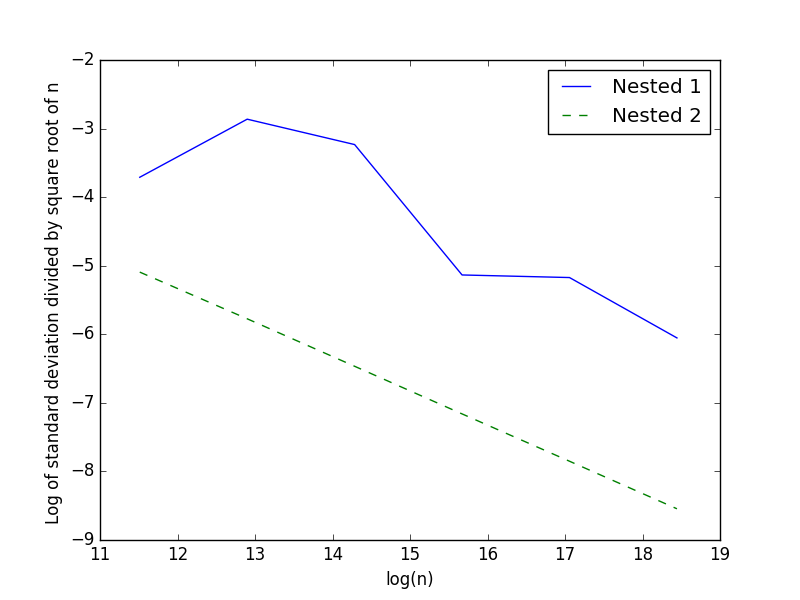}
  \caption{Estimation and error in $d=4$  on case test A. Maturity $T=1.5$}
\label{semiFig2}
\end{figure}

\begin{figure}[H]
  \centering
  \includegraphics[width=7cm]{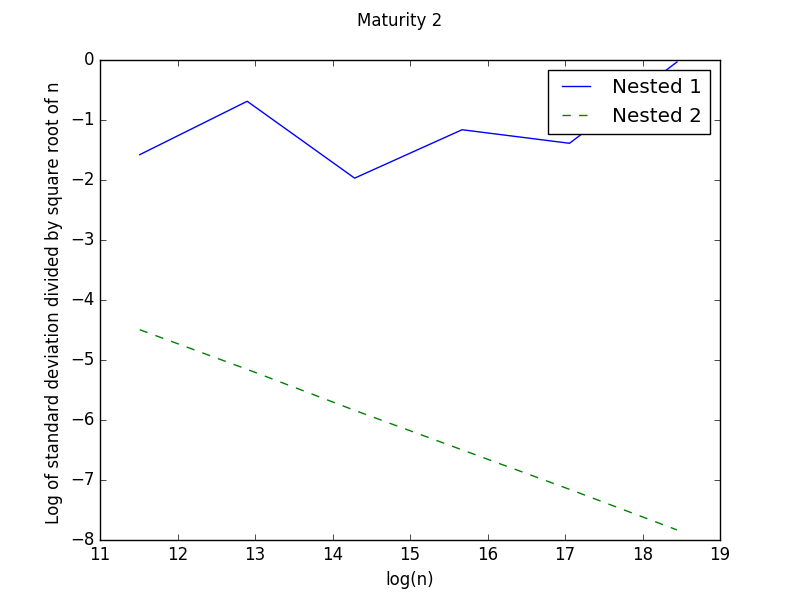}
  \includegraphics[width=7cm]{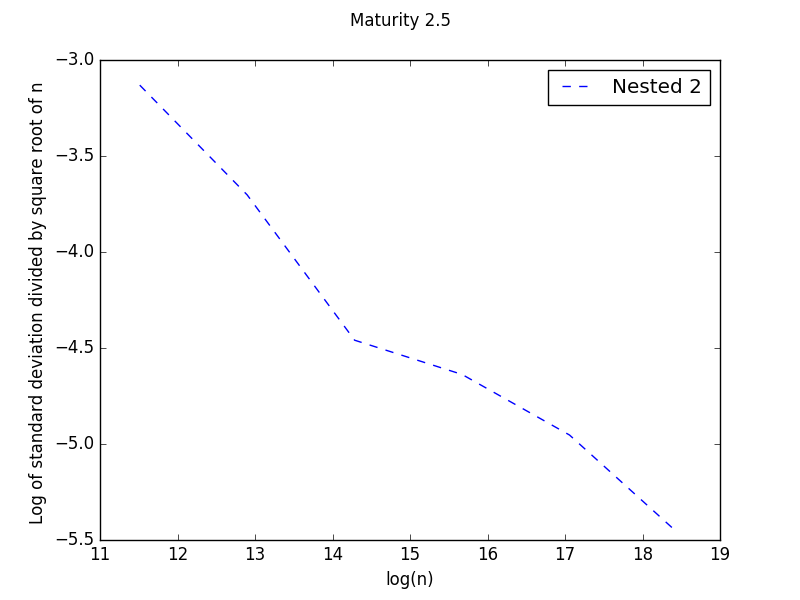}
  \caption{Error in $d=4$  on case test A. Maturity $T=2.$, $T=2.5$}
\label{semiFig3}
\end{figure}

We first give on figures \ref{semiFig1},  \ref{semiFig2} and  \ref{semiFig3}  the results obtained for test A for different maturities and a dimension $d=4$ so the analytical solution is $-0.508283$. We plot for each maturity :
\begin{itemize}
\item the solution obtained by increasing the number of Monte Carlo scenarios used,
\item the error calculated as explained in the introduction.
\end{itemize}
Nested $n$  curves stand for the curves using  the nested method of order  $n$, so the  Nested $1$ curve stands for the original method.
On figure \ref{semiFig3}, for maturity $2.5$  the error observed with the orignal method (Nested 1) is around 1000 so it has not been plotted.
Because of the number of branching due to the gamma law, it seems difficult to use a nested method of order $n>2$ for long maturities : the time needed explodes.
But clearly the nested method permits to have accurate solution for longer maturities.

\begin{figure}[H]
  \centering
  \includegraphics[width=7cm]{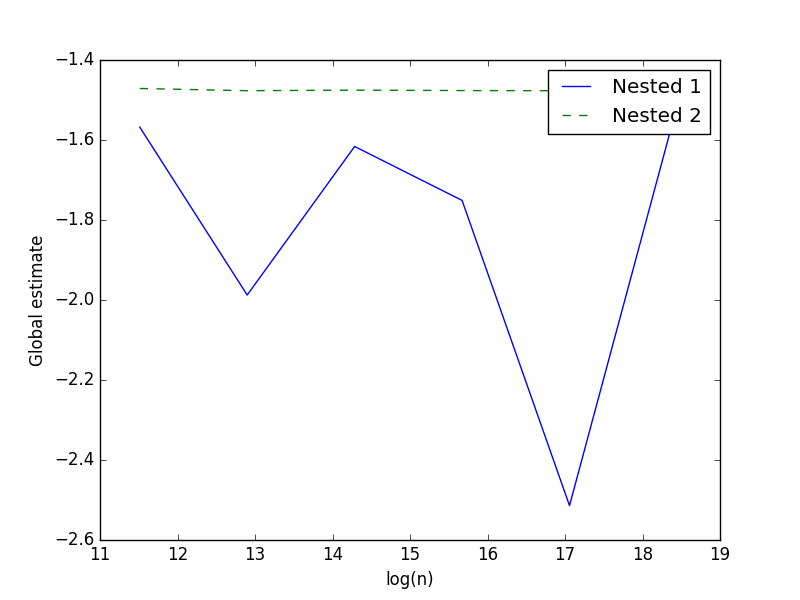}
 \includegraphics[width=7cm]{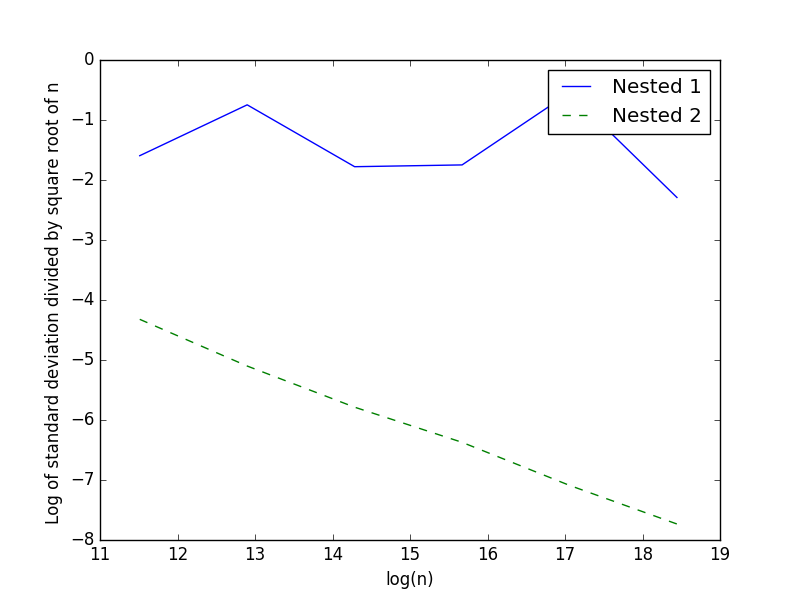}
   \caption{Estimation and error in $d=6$  on case test A. Maturity $T=2.$}
\label{semiFig3_1}
\end{figure}

For a maturity of $2$ we also give the results obtained in dimension $6$ on figure \ref{semiFig3_1} giving an analytical solution  $-1.4769$ : once again the original method fails to converge while the nested one give good results.

\subsection{Adaptation of the original branching to the re-normalization technique}
\label{sec:SemiLinRenorm}
As introduced in \cite{LTTW}, we introduce a modification of the original branching process that let us use exponential laws  for the branching dates 
to treat the $Du$ terms in the method previously described.
 	Recall that $\Kcb^1_T = \{(1)\}$, we introduce an associated ghost particle, denoted by $(1^{1})$,
	and denote $ \Kct^1_T := \{(1), (1^{1}) \}$.
	Next, given the collection $\Kct^n_T$ of all particles (as well as ghost particles) of generation $n$,
	we define the collection $\Kct^{n+1}_T$ as follows.
	For every $k = (k_1, \cdots, k_n) \in \Kct^n$, we denote by $o(k) = (\hat k_1, \cdots, \hat k_n)$ its original particles,
	where $\hat k_i := j$ when $k_i = j ~\mbox{or}~ j^{1}$.
	Further, when $k = (k_1, \cdots, k_n)$ is such that $k_n \in \N$, we denote $k^{1} := (k_1, \cdots, k_{n-1}, k_n^{1})$.
	The mark of $k \in \Kct^n$ will be the same as its original particle $o(k)$, i.e. $\theta_k := \theta_{o(k)}$;
	and $T_k := T_{o(k)}$, $\Delta T_k := \Delta T_{o(k)}$ and $\tau^k = \tau^{o(k)}$.
	Define also $\Kch^n_T:= \{k \in \Kct^n_T ~: o(k) \in \Kc_T\}$.	
	For every $k = (k_1, \cdots, k_n) \in \Kct^n_T \setminus \Kch^n_T$, we still define the set of its offspring particles by
	$$S(k) := \{(k_1, \cdots, k_n, 1), \cdots, (k_1, \cdots, k_n, L) \},$$
	and the set of ghost offspring particles by
	$$
		S^{1}(k)
		~~:=~~
		\big\{ (k_1, \cdots, k_n, 1^{1}), \cdots, (k_1, \cdots, k_n, L^{1}) \big\}.
	$$
 	Then the collection $ \Kct^{n+1}_T$ of all particles (and ghost particles) of generation $n+1$ is
	$$
		\Kct^{n+1}_T
		~:=~
		\cup_{k \in \Kct^n_T \setminus \Kch^n_T} \big( S(k) \cup S^{1}(k) \big).
	$$
	Define also
	$$
		\Kct_T := \cup_{n \ge 1} \Kct^n_T,
		~~~\mbox{and}~
		\Kch_T := \cup_{n \ge 1} \Kch^n_T.
	$$

\subsubsection{The original re-normalization technique}
\label{subsec:renorm}
        We next equip each particle with a Brownian motion in order to define a branching Brownian motion.
	Further, let $W^{\emptyset}_0 := 0$, and for every $k = (k_1, \cdots, k_n) \in \Kct^n_T$, let
	\begin{flalign}\label{eq:brownRenorm}
		W^{k}_s
		~:=~
		W^{k-}_{T_{k-}}
		~+~
		\1_{k_n \in \N}
		\hat W^{o(k)}_{s - T_{k-}},
		~~~\mbox{and}~~
		X^{k}_s := \mu s +\sigma_0  W^k_s,
		~~~\forall s \in [T_{k-}, T_k].
	\end{flalign}

\begin{figure}[H]
\centering
 \begin{subfigure}[b]{0.45\textwidth}
   \includegraphics[width=\textwidth]{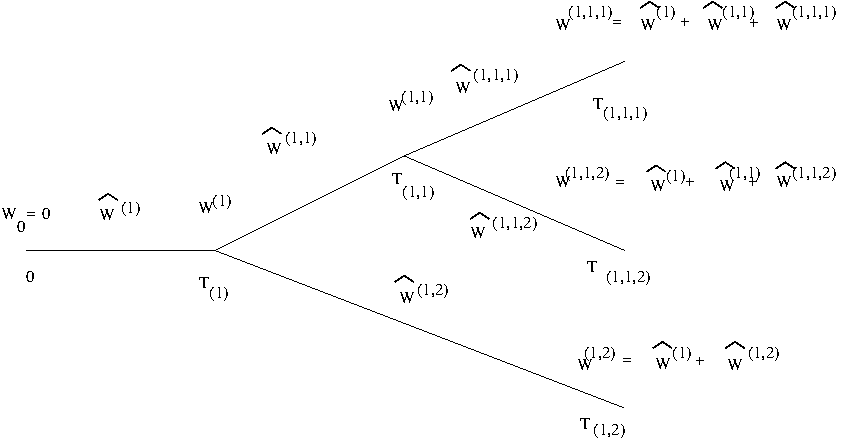}
   \caption{Original Galton-Watson tree \\
	{\tiny
	\begin{tabular}{llcll}
	\hline
	$W^{(1)} = \hat W^{(1)}$ \\ \hline
	$W^{(1,1)} = \hat W^{(1)} + \hat W^{(1,1)}$ \\ \hline
	$W^{(1,2)} = \hat W^{(1)} + \hat W^{(1,2)}$ \\ \hline
	$W^{(1,1,1)} = \hat W^{(1)} + \hat W^{(1,1)}+ \hat W^{(1,1,1)}$ \\ \hline
	$W^{(1,1,2)} = \hat W^{(1)} + \hat W^{(1,1)}+ \hat W^{(1,1,2)}$ \\
	\hline
	\end{tabular}}
}
 \end{subfigure}
 \begin{subfigure}[b]{0.45\textwidth}
   \includegraphics[width=\textwidth]{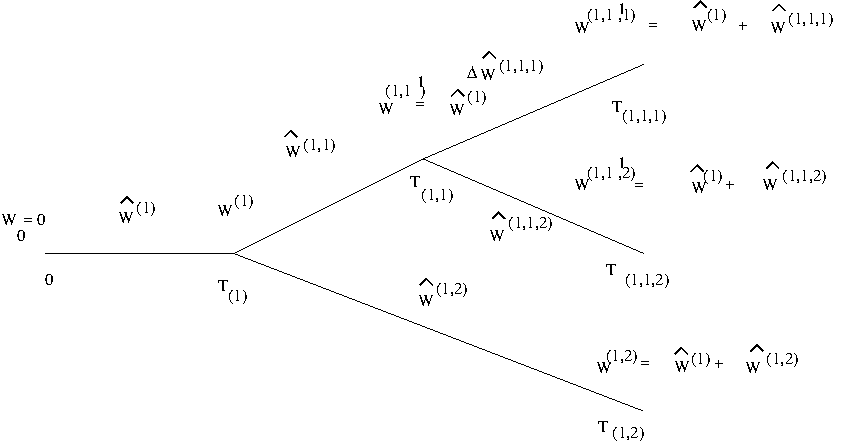}
   \caption{Tree with ghost particle $k=(1,1^{1})$ \\
	{\tiny
	\begin{tabular}{llcll}
	\hline
	$W^{(1)} = \hat W^{(1)}$ \\ \hline
	$W^{(1,1^{1})} = \hat W^{(1)} $ \\ \hline
	$W^{(1,2)} = \hat W^{(1)} + \hat W^{(1,2)}$ \\ \hline
	$W^{(1,1^{1},1)} = \hat W^{(1)} + \hat W^{(1,1,1)}$ \\ \hline
	$W^{(1,1^{1},2)} = \hat W^{(1)} + \hat W^{(1,1,2)}$ \\
	\hline
	\end{tabular}}
   }
 \end{subfigure}
 \begin{subfigure}[b]{0.45\textwidth}
\includegraphics[width=\textwidth]{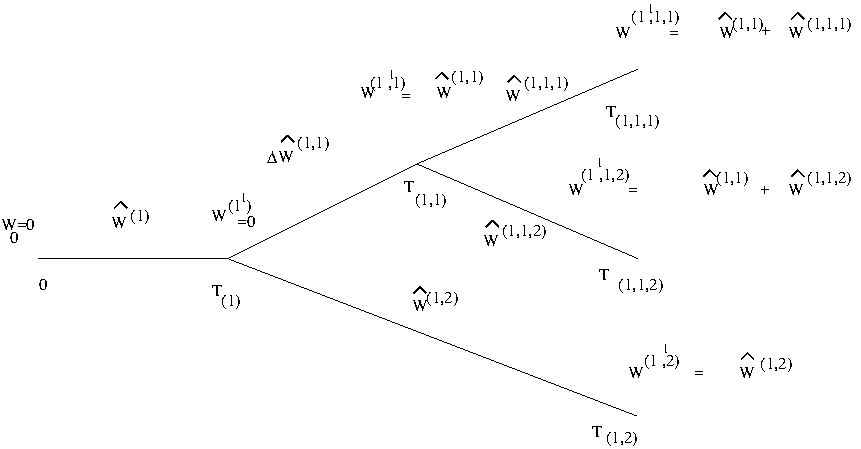}
\caption{Tree with ghost particle $k=(1^{1})$\\
	{\tiny
	\begin{tabular}{llcll}
	\hline
	$W^{(1^{1})} = 0$ \\ \hline
	$W^{(1^{1},1)} = \hat W^{(1,1)}$ \\ \hline
	$W^{(1^{1},2)} = \hat W^{(1,2)}$ \\ \hline
	$W^{(1^{1},1,1)} = \hat W^{(1,1)} + \hat W^{(1,1,1)}$ \\ \hline
	$W^{(1^{1},1,2)} = \hat W^{(1,1)} + \hat W^{(1,1,2)}$ \\
	\hline
	\end{tabular}}
}
\end{subfigure}
\begin{subfigure}[b]{0.45\textwidth}
\includegraphics[width=\textwidth]{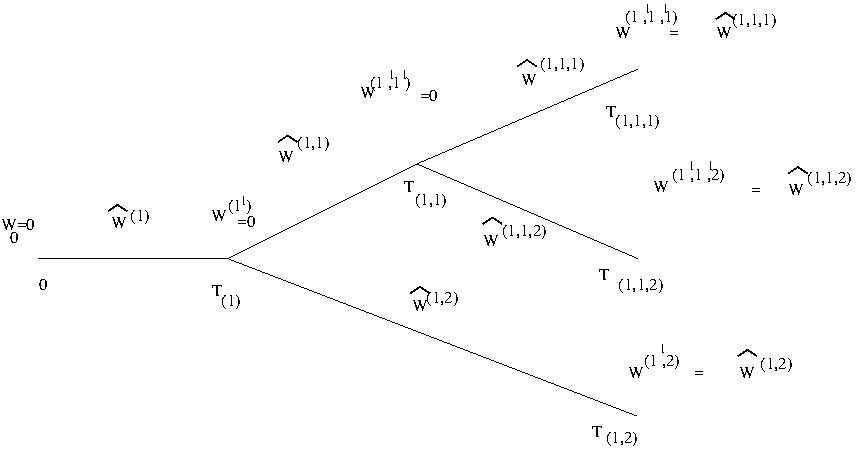}
\caption{Tree with ghost particles $k=(1^{1})$ and $k=(1^{1},1^{1})$\\
	{\tiny
	\begin{tabular}{llcll}
	\hline
	$W^{(1^{1})} = 0$ \\ \hline
	$W^{(1^{1},1^{1})} = 0$ \\ \hline
	$W^{(1^{1},2)} = \hat W^{(1,2)}$ \\ \hline
	$W^{(1^{1},1^{1},1)} = \hat W^{(1,1,1)}$ \\ \hline
	$W^{(1^{1},1^{1},2)} = \hat W^{(1,1,2)}$ \\
	\hline
	\end{tabular}}
}
\end{subfigure}
\caption{Original Galton-Watson tree, different trees with ghost particles 
	(excluding ghost particles at the extreme leaves) for a Brownian motion where
	$W^k$ stands for $W^k_{T_k}$ and  $\hat W^k$ stands for $\hat W^k_{\Delta T_k}$.
}
\label{figTree}
\end{figure}

On figure \ref{figTree}, we give the original Galton-Watson tree and the ghost particles associated.

The initial equation \eqref{eq:uVal} remains unchanged (first  step of the algorithm) but equation \eqref{eq:RecVal} is modified
by replacing the term 
\begin{flalign*}
\E_{T_{(1)},X_{T_{(1)}}}\big[ (\sigma_0^\top)^{-1} \frac{\hat W^{(1,p)}_{\Delta T_{(1,p)}}}{\Delta T_{(1,p)}} \phi\big(T_{(1,p)},X^{(1,p)}_{T_{(1,p)}}\big)  \big]
\end{flalign*}
by
\begin{flalign}
\label{eq:folRep}
 \E_{T_{(1)},X_{T_{(1)}}}\big[ (\sigma_0^\top)^{-1} \frac{\hat W^{(1,p)}_{\Delta T_{(1,p)}}}{\Delta T_{(1,p)}} \big(\phi\big(T_{(1,p)},X^{(1,p)}_{T_{(1,p)}}\big)  -\phi\big(T_{(1,p)},X^{(1,p^1)}_{T_{(1,p)}}\big) \big)  \big].
\end{flalign}
Notice that since $W^{(1,p^1)}$ has been obtained  by \eqref{eq:brownRenorm}, ${\hat W^{(1,p)}}_{\Delta T_{(1,p)}}$ and $\phi\big(T_{(1,p)},W^{(1,p^1)}_{T_{(1,p)}}\big)$ are orthogonal so that adding the second term acts as a control variate.
Recursively using the modified version of equation \eqref{eq:RecVal} induced by the use of \eqref{eq:folRep}, \cite{LTTW} gave defined the re-normalized estimator by a backward induction:
let $\psih_k := \frac{g(X^k_T)}{\Fb(\Delta T_k)}$ for every $k \in \Kch_T$, then let
\begin{flalign}\label{eq:ghostRep}
\psih_k
:=
\frac{1}{\rho(\Delta T_k) } \big( h(T_k, X^k_{T_k})+c(T_k, X^k_{T_k})
\prod_{{\tilde k} \in S(k)} \!\!\! \big(\psih_{{\tilde k}} -  \psih_{{\tilde k}^{1}} \1_{\{\theta({\tilde k}) \neq 0\}}\big) \Wc_{{\tilde k}} \big),
		~~\mbox{for}~k \in \Kct_T \setminus \Kch_T.
\end{flalign}
where the weights $ $ are given by equation \eqref{eq:weigthSem},
so we have 
\begin{flalign*}
u(0,x)
 = \E_{0,x}\Big[  \psih_{(1)}\Big].
\end{flalign*}

As explained in section \ref{sec:origAlgo}, equation \eqref{eq:RecVal} used in representation \eqref{eq:initRep} force us to take laws for branching dates  with a high probability  of low values that leads to a high number of recursions defined by equation  \eqref{eq:backRep}.
Besides such laws using some rejection algorithm, as gamma laws, are very costly to generate.
The use of  \eqref{eq:folRep} permits us to use exponential laws very cheap to simulate and with a low probability of small values. \\
Indeed it can be easily seen in the linear case  ($f$  function of the gradient with $\ell_0=0$, $m=1$ and $\ell_1=1$) by conditioning with respect to the number of branching that the variance is bounded for small maturities
and coefficients if
\begin{flalign}
\label{condRenorm}
\E_{0,x}\big[ \big(\psih_{k} -  \psih_{k^{1}} \1_{\{\theta(k) \neq 0\}}\big)^2 \frac{ \left(b_{\theta_k}(T_{k-}, X^k_{T_{k-}})
		\cdot (\sigma_0^{\top})^{-1} \hat W^{o(k)}_{\Delta T_k}\right)^2}{(\Delta T_k)^2} \big] < \infty.
\end{flalign}
By $X^{k^{1}}_t$ construction using $g$ regularity, it is easily seen that for small time steps $\Delta T_k$,  $\E_{0,x,\Delta T_k}\big[ (\psih_{k} -  \psih_{k^{1}})^2 \big] = O(\Delta T_k)$ as $\Delta T_k \longrightarrow 0$ and \eqref{condRenorm} is satisfied for every $\rho$ densities.\\

\subsubsection{Re-normalization techniques and antithetic}
\label{subsec:renormAnt}
We give a version of the re-normalization technique using antithetic variables.
Equation \eqref{eq:brownRenorm} is modified  by :
\begin{flalign}\label{eq:brownRenormAnti}
W^{k}_s
~:=~
W^{k-}_{T_{k-}}
		~+~
		\1_{k_n \in \N}
		\hat W^{o(k)}_{s - T_{k-}}  ~-~ \1_{k_n \notin \N}
		\hat W^{o(k)}_{s - T_{k-}}, \nonumber \\
		~~~\mbox{and}~~
		X^{k}_s := \mu s +\sigma_0  W^k_s,
		~~~\forall s \in [T_{k-}, T_k],
                \end{flalign}
 for every $k = (k_1, \cdots, k_n) \in \Kct^n_T$.\\
Then equation \eqref{eq:RecVal} is modified by  :
\begin{itemize}
\item First , replacing the term tacking into account the power of $u$
\begin{flalign*}
\phi\big(T_{(1,j)},X^{(1,j)}_{T_{(1,j)}}\big)
\end{flalign*}
 by
\begin{flalign*}
 \frac{1}{2} \big( \phi\big(T_{(1,j)},X^{(1,j)}_{T_{(1,j)}}\big)+
\phi\big(T_{(1,j)},X^{(1,j^1)}_{T_{(1,j)}}\big) \big),
\end{flalign*}
\item and the term taking into account the gradient 
\begin{flalign*}
 \E_{T_{(1)},X_{T_{(1)}}}\big[ (\sigma_0^\top)^{-1} \frac{\hat W^{(1,p)}_{\Delta T_{(1,p)}}}{\Delta T_{(1,p)}} \phi\big(T_{(1,p)},X^{(1,p)}_{T_{(1,p)}}\big)  \big]
\end{flalign*}
by 
\begin{flalign*}
\E_{T_{(1)},X_{T_{(1)}}}\big[ (\sigma_0^\top)^{-1} \frac{\hat W^{(1,p)}_{\Delta T_{(1,p)}}}{\Delta T_{(1,p)}} \frac{1}{2}\big(\phi\big(T_{(1,p)},X^{(1,p)}_{T_{(1,p)}}\big)  -\phi\big(T_{(1,p)},X^{(1,p^1)}_{T_{(1,p)}}\big) \big)  \big].
\end{flalign*} 
\end{itemize}

Notice that with this version the variance of the gradient term is finite with the same argument as in the original re-normalization version in subsection \ref{subsec:renorm}.\\
By backward induction we get the re-normalized antithetic estimator  modifying \eqref{eq:ghostRep} by:
\begin{flalign}\label{eq:ghostRepAntithetic}
\psih_k
& := 
\frac{1}{\rho(\Delta T_k) } \big( h(T_k, X^k_{T_k})+c(T_k, X^k_{T_k})
\prod_{{\tilde k} \in S(k)} \!\!\! \frac{1}{2} \big(\psih_{{\tilde k}} -  \psih_{{\tilde k}^{1}} \1_{\{\theta({\tilde k}) \neq 0\}} + \psih_{{\tilde k}^{1}} \1_{\{\theta({\tilde k}) = 0\}} \big) \Wc_{{\tilde k}} \big),
	\nonumber \\
&  \mbox{for}~k \in \Kct_T \setminus \Kch_T.
\end{flalign}
where the weights $ $ are given by equation \eqref{eq:weigthSem}.
Then  we have 
\begin{flalign*}
u(0,x)
 = \E_{0,x}\Big[  \psih_{(1)}\Big].\end{flalign*}

\subsubsection{Numerical result for semi linear with re-normalization}

We apply our nested algorithm on the original re-normalized technique  and on the re-normalization technique with antithetic variables on  two test cases.\\
First we give some results for test case A in dimension 4.
We give the Monte Carlo error  obtained by the nested method on figure \ref{semiFig4}.
For the maturity $T=3$, without nesting the error of the original re-normalization technique  has an order of magnitude of 2000 so the curve has not been given.
For the maturity $T=4$, the nested  original re-normalization technique with  an order 2 doesn't seem to converge.
\begin{figure}[H]
  \centering
  \includegraphics[width=7cm]{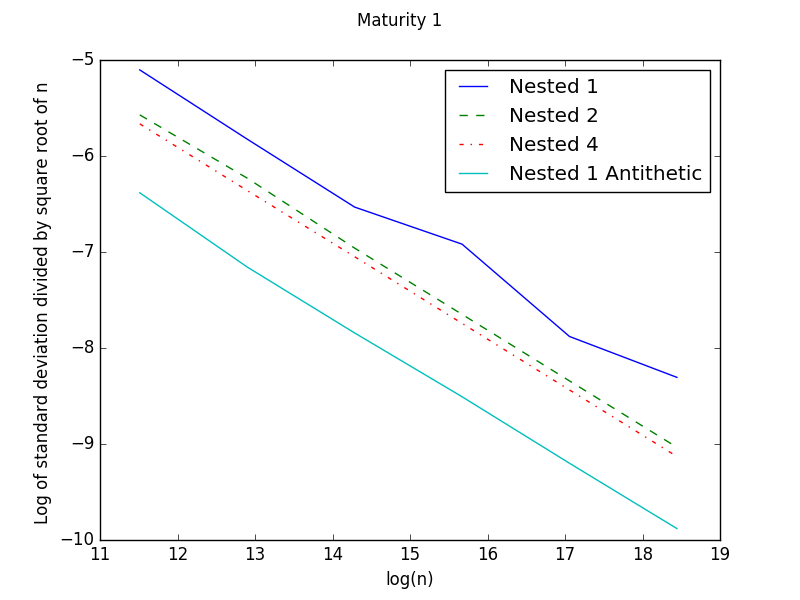}
  \includegraphics[width=7cm]{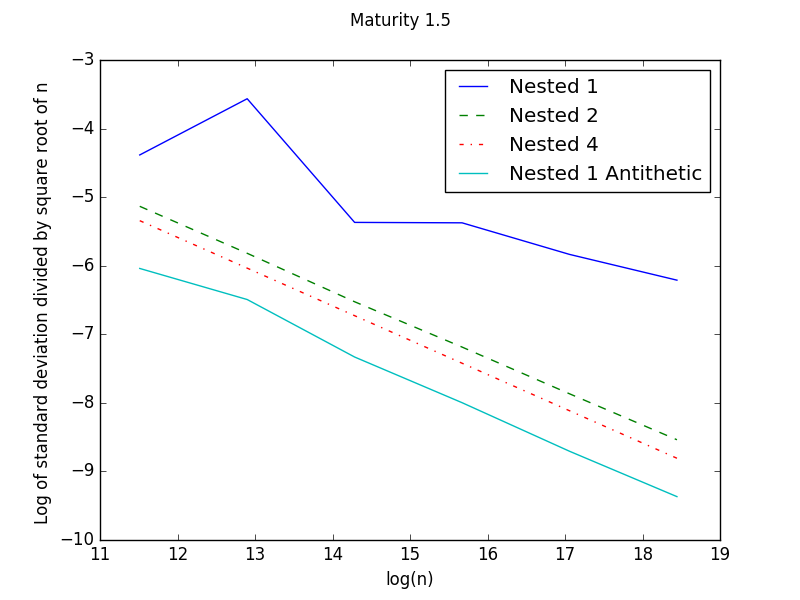}
  \includegraphics[width=7cm]{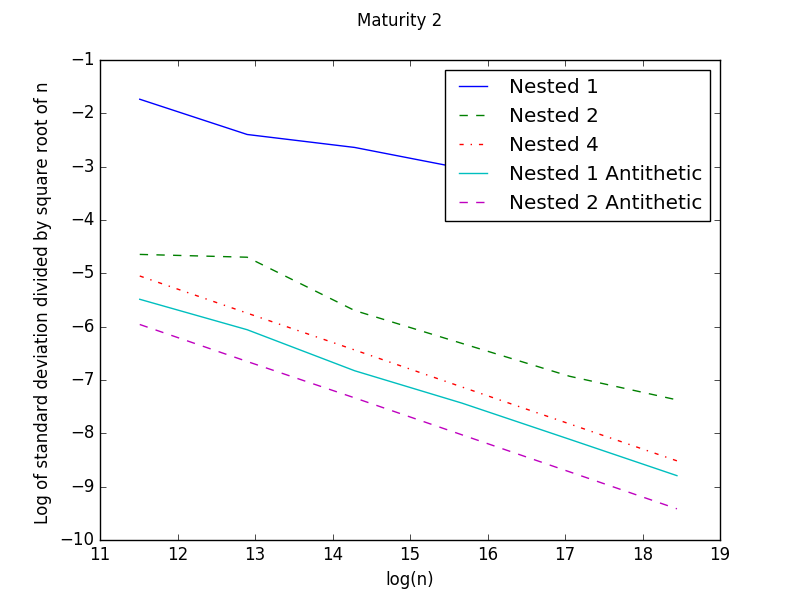}
  \includegraphics[width=7cm]{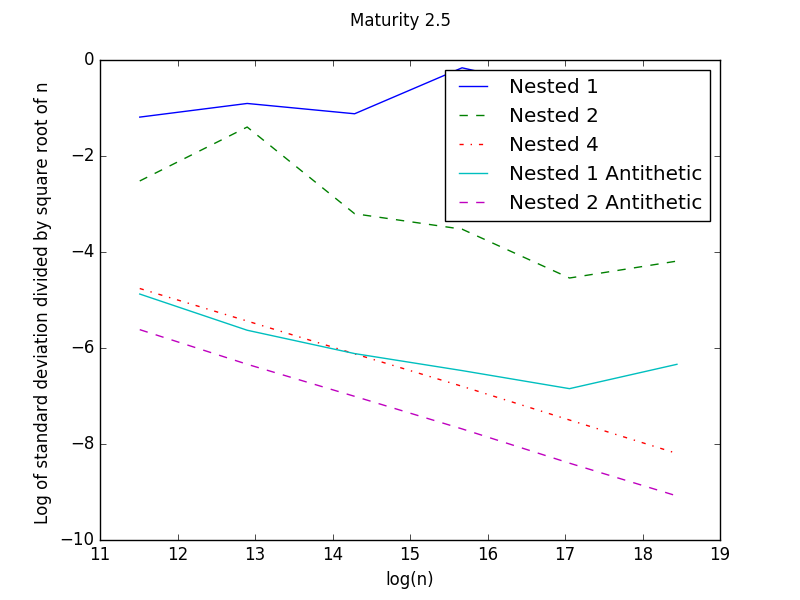}
  \includegraphics[width=7cm]{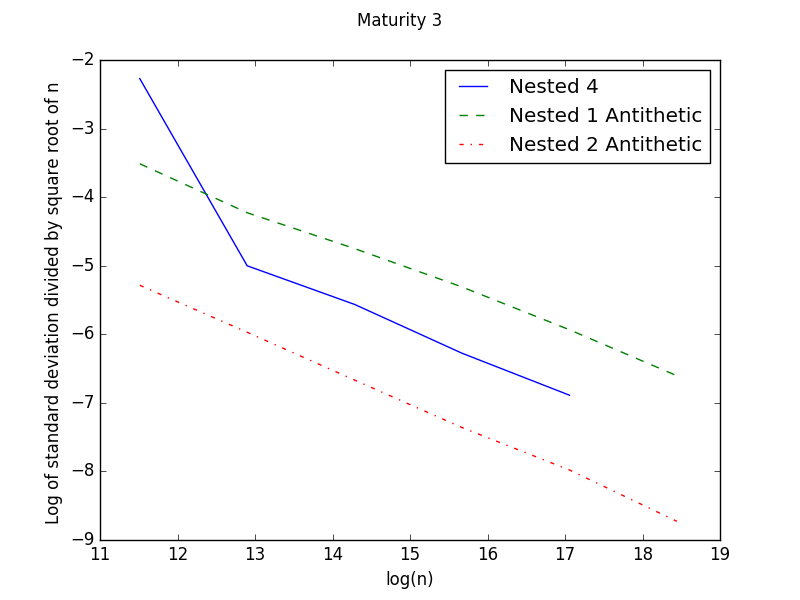}
  \includegraphics[width=7cm]{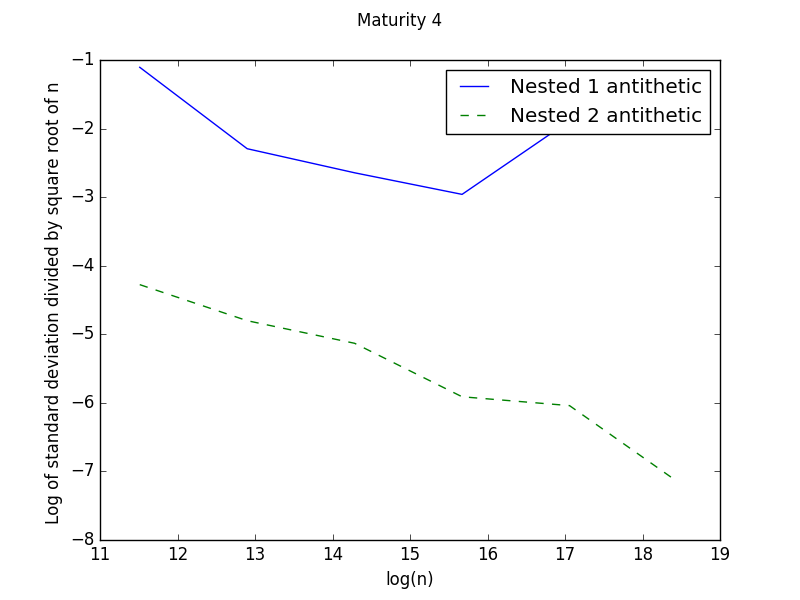}
   \caption{Error in $d=4$  on case test A for different maturities}
\label{semiFig4}
\end{figure}
As the maturity increases, nesting with a higher order becomes necessary.
Notice that with the re-normalization it is possible to use the nested method of a high order because of the small number of branching used.
For example, for $T=2$,  for an accuracy of $0.0004$, in dimension $d=4$:
\begin{itemize}
\item the original method in section \ref{sec:origAlgo} with a nested method of order 2 achieves an accuracy of $0.0004$ for a CPU time of $1500$ seconds using 28 cores,
\item the re-normalized version of section \ref{subsec:renorm} with a nested method of order 4 reaches the same accuracy in $1800$ seconds,
\item the re-normalized version with antithetic of section \ref{subsec:renormAnt} without nesting  reaches the same accuracy in $11$ seconds.
\end{itemize}

\begin{figure}[H]
  \centering
\includegraphics[width=7cm]{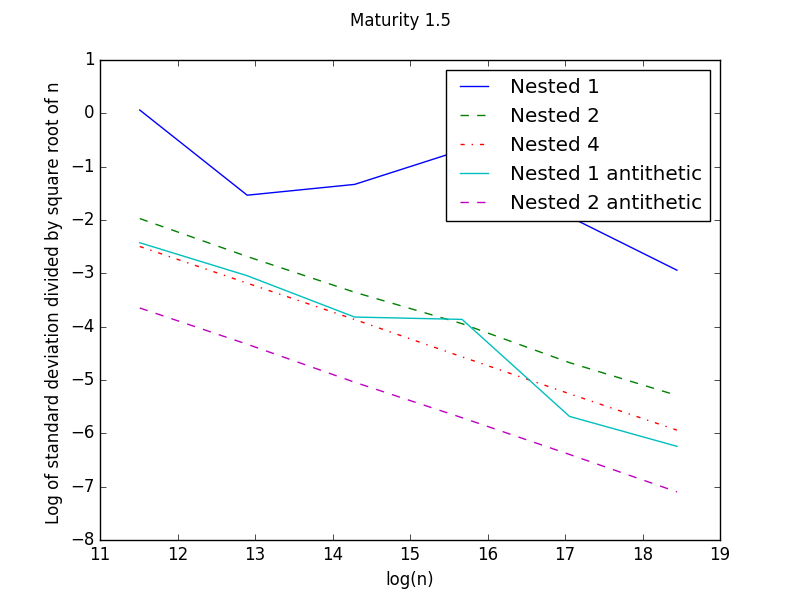}
\caption{Error in $d=6$ for the term  $b.Du$  on  case test A for $T=1.5$.}
\label{semiFig4_}
\end{figure}

\begin{figure}[H]
  \centering
  \includegraphics[width=7cm]{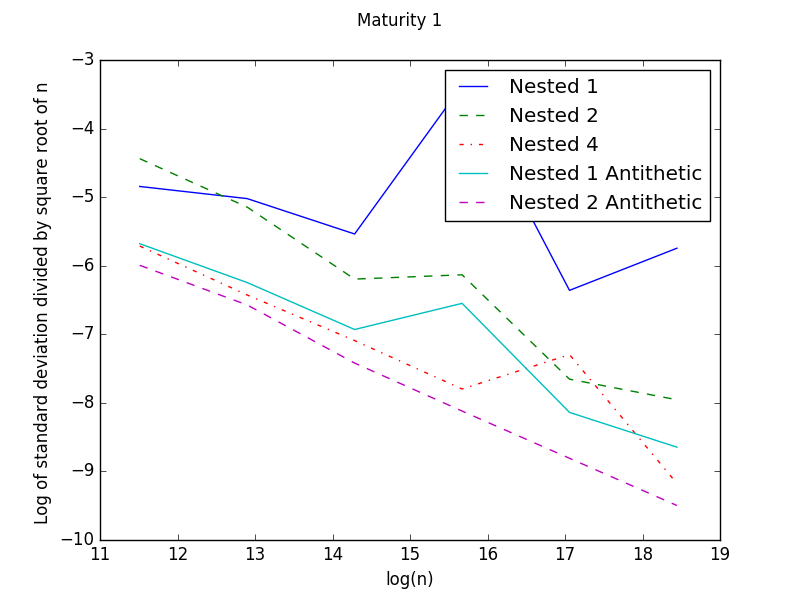}
  \includegraphics[width=7cm]{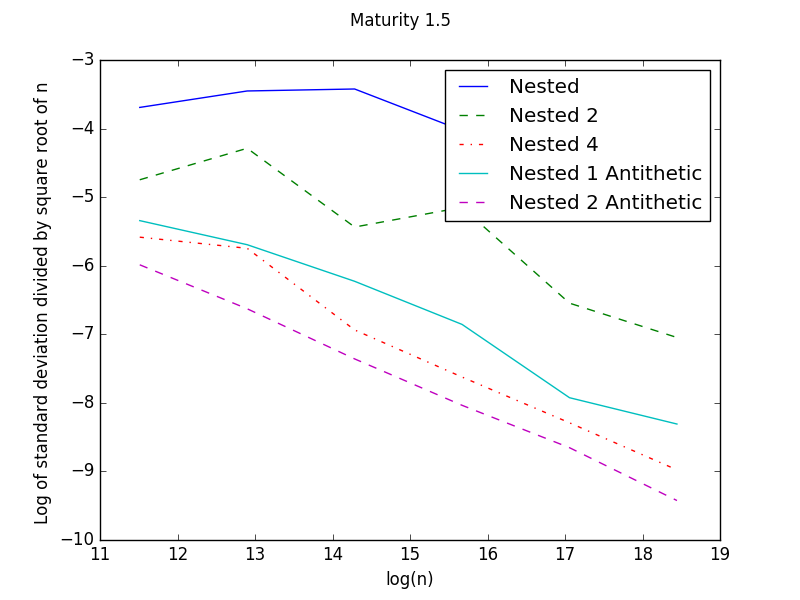}
  \includegraphics[width=7cm]{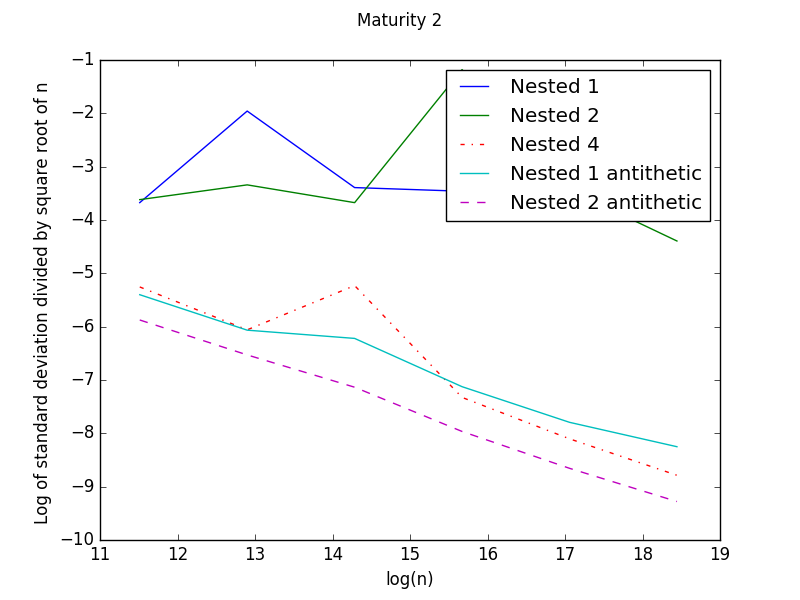}
  \includegraphics[width=7cm]{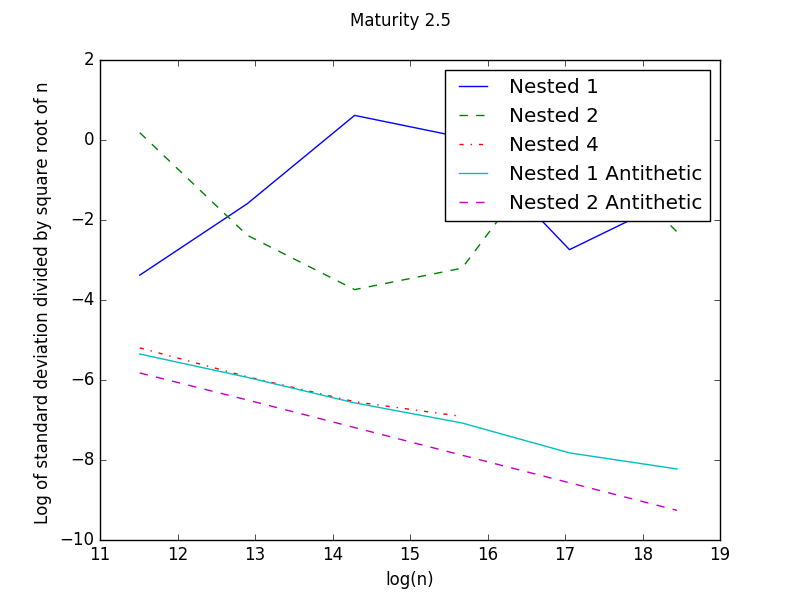}
   \caption{Error in $d=4$  on case test B for different maturities.}
\label{semiFig5}
\end{figure}

\begin{figure}[H]
\centering
  \includegraphics[width=7cm]{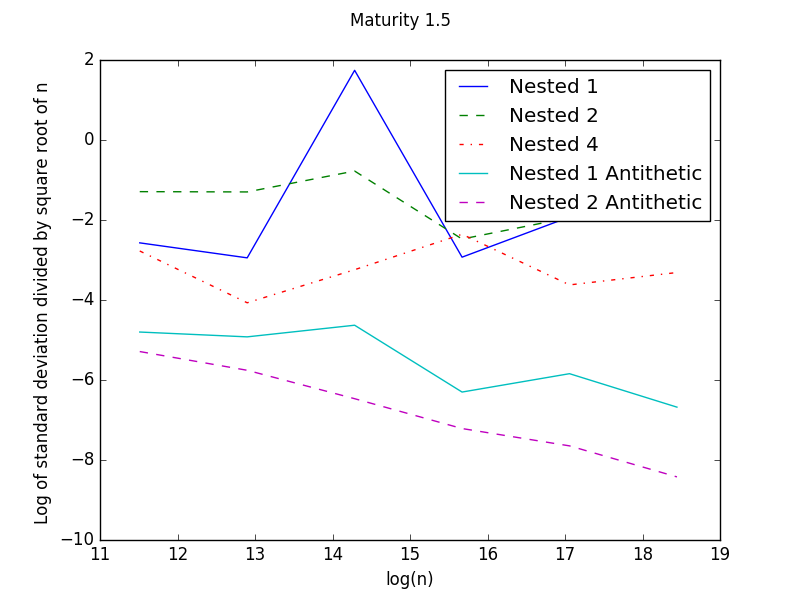}
  \includegraphics[width=7cm]{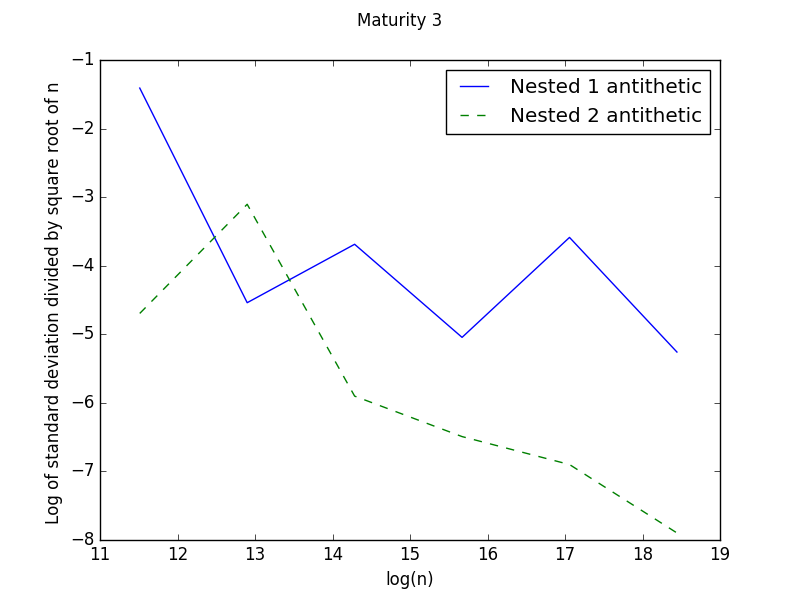}
   \caption{Error in $d=6$  on case test B  for different maturities.}
\label{semiFig6}
\end{figure}

For the same test case A we plot in dimension 6 the error on figure \ref{semiFig5} to show that the method converges in high dimension.
\begin{figure}[H]
  \centering
  \includegraphics[width=7cm]{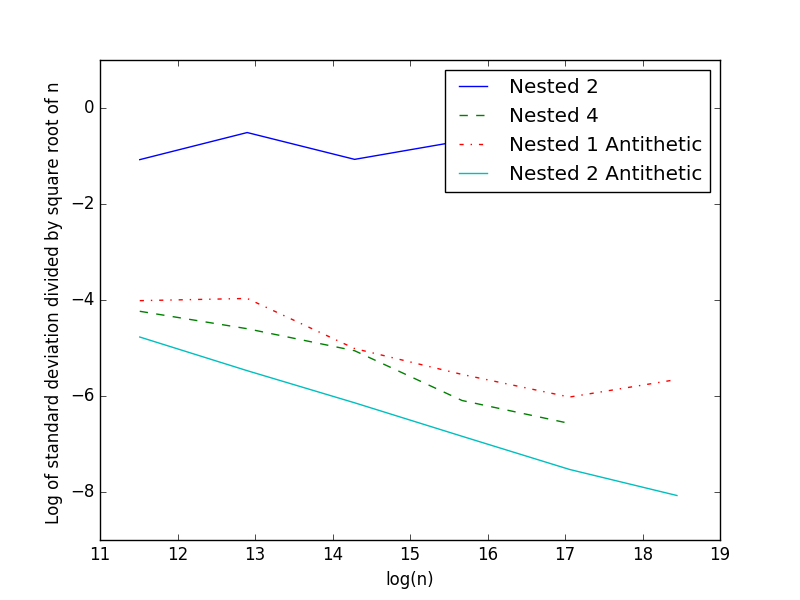}
  \caption{Error in $d=6$  on  case test A for $T=3$.}
\label{semiFig4}
\end{figure}
Besides on figure \ref{semiFig4_}, we show that the derivative is accurately calculated.
We then use a second test case B :
For a given dimension $d$ , we take $\sigma_0 = \frac{1}{\sqrt{d}} \I_d$, $\mu= \0$,
\begin{flalign*}
 f(t,x,y,z)=  \frac{0.1}{d}  (\1 \cdot z)^2
\end{flalign*}
 with a terminal condition $g(x) = \cos( x_1+ \cdots + x_d)$.
This test case cannot be solve by the nested method without re-normalization due to the high cost involved by the potential high number of branching. We give the results obtained for case B by the re-normalization methods of section  \ref{subsec:renorm} and \ref{subsec:renormAnt}  in dimension 4 on figure \ref{semiFig5}.
At last we give the results obtained in dimension $6$ pour $T=1.5$ and $T=3$ on figure \ref{semiFig6}.
The nested method with re-normalization and antithetic appears to be the most effective and permits to solve semi-linear equations with quite long maturities.
The re-normalization technique is however far more memory consuming than the original scheme of section \ref{sec:origAlgo}. This memory cost  explodes with very high maturities.
The nested version of the   original scheme of section \ref{sec:origAlgo} isn't affected by these  memory problems but is affected by an explosion of the computational time with longer maturities.
\subsection{Extension to variable coefficients}
\begin{figure}[H]
  \centering
  \includegraphics[width=7cm]{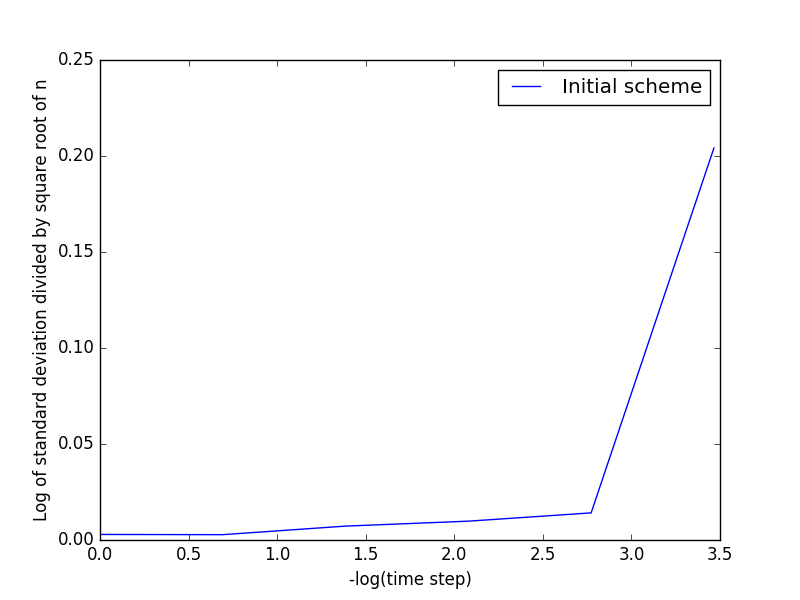}
  \includegraphics[width=7cm]{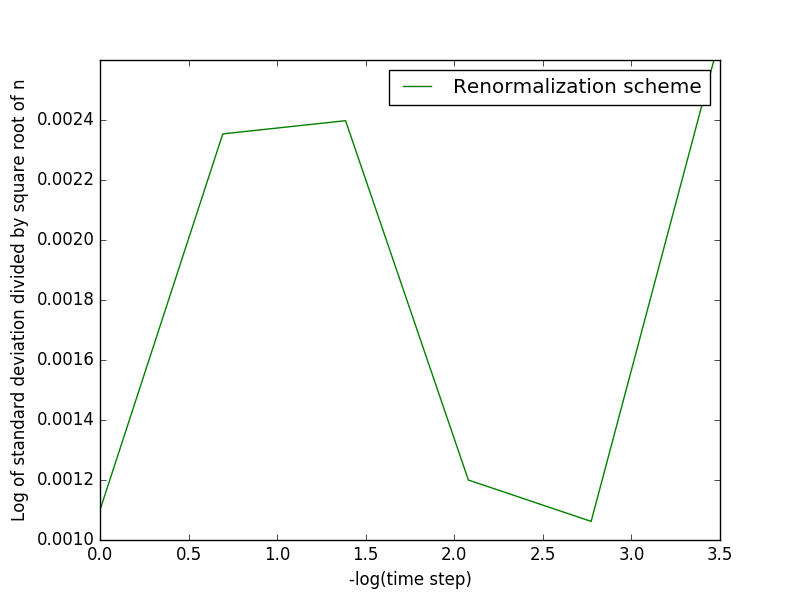}
\caption{A Burgers case in dimension 4 : comparison of Euler schemes error for  the original method and the renormalization method.}
\label{figEuler}
\end{figure}

In the case of time and space dependent  coefficients $\mu$ and $\sigma_0$ of the PDE, it is possible to use the method consisting in ``freezing'' the coefficients first proposed  in \cite{HTT2} for non fixed $\mu$  and extended in  the general case in \cite{DOW}.
This method increases the variance of the estimator, therefore it is more efficient for treating log maturities to use an Euler scheme to take into account the variation of the coefficients.
Introducing an Euler time step $\delta t$, between the dates $T_{k-}$ and  $T_k$, the SDE is discretized as :
\begin{align*}
X^k_{T_{k-}+i \delta t }=& X^k_{T_{k-}+(i-1) \delta t} + \mu(T_{k-}+(i-1) \delta t , X^k_{T_{k-}+(i-1) \delta t}) \delta t + \\
                   &  \sigma_0(T_{k-}+(i-1) \delta t, X^k_{T_{k-}+(i-1) \delta t}) \hat  W^{k,i}_{\delta t},
\mbox{ for } i=1, ..,N,\\
X^k_{T_k} =& X^k_{T_{k-}+ N \delta t}+\mu(T_{k-}+N \delta t, X^k_{T_{k-}+ N \delta t}) (\Delta T_k - N\delta t) +  \\
             &      \sigma_0(T_{k-}+N \delta t, X^k_{T_{k-}+N \delta t}) \hat  W^{k,i}_{\Delta T_k - N\delta t},
\end{align*}

where $N = \lfloor \frac{\Delta T_k}{\delta t} \rfloor$,  and 
$(\hat W^{k,i})_{k = (k_1, \cdots, k_{n-1}, k_n) \in \N^n, n>1, i\ge 1}$  is a sequence of independent $d$-dimensional Brownian motion.\\
Using an integration by part on the first time step, in the original scheme of section \ref{origScheme}, the gradient term in equation \eqref{eq:RecVal} is replaced
\begin{align}
\label{eq:eulerOrig}
\E_{T_{(1)},X_{T_{(1)}}}\big[ (\sigma_0(T_{(1)},X^{(1)}_{T_{(1)}})^\top)^{-1}\frac{ \hat W^{(1,p),1}_{\min(\delta t,\Delta T_{(1,p)})}}{\min(\delta t,\Delta T_{(1,p)})} \phi\big(T_{(1,p)},X^{(1,p)}_{T_{(1,p)}}\big)\big]
\end{align}
In the case of the renormalization technique of section \ref{subsec:renorm}, the ghost is obtained from the original particule by removing the part associated to the first brownian.
Then for every $k = (k_1, \cdots, k_n) \in \Kct^n_T$, the particule dynamic is given by  
\begin{align*}
X^k_{T_{k-}+  \delta t }: =& X^k_{T_{k-}} + \mu(T_{k-}, X^k_{T_{k-}}) \delta t + \\
                   &  \1_{k_n \in \N} \sigma_0(T_{k-}, X^k_{T_{k-}}) \hat  W^{k,1}_{\delta t},\\
X^k_{T_{k-}+i \delta t }=& X^k_{T_{k-}+(i-1) \delta t} + \mu(T_{k-}+(i-1) \delta t , X^k_{T_{k-}+(i-1) \delta t}) \delta t + \\
                   &  \sigma_0(T_{k-}+(i-1) \delta t, X^k_{T_{k-}+(i-1) \delta t}) \hat  W^{k,i}_{\delta t},
\mbox{ for } i=2, ..,N,\\
X^k_{T_k} =& X^k_{T_{k-}+ N \delta t}+\mu(T_{k-}+N \delta t, X^k_{T_{k-}+ N \delta t}) (\Delta T_k - N\delta t) +  \\
             &      \sigma_0(T_{k-}+N \delta t, X^k_{T_{k-}+N \delta t}) \hat  W^{k,i}_{\Delta T_k - N\delta t},
\end{align*}
if $N>0$ and
\begin{align*}
X^k_{T_{k}}: =& X^k_{T_{k-}} + \mu(T_{k-}, X^k_{T_{k-}}) \Delta T_k + \\
                   &  \1_{k_n \in \N} \sigma_0(T_{k-}, X^k_{T_{k-}}) \hat  W^{k,1}_{\Delta T_k}
\end{align*}
otherwise. \\
The renormalization technique of section \ref{subsec:renorm} leads to the following estimation of the gradient  in equation \eqref{eq:RecVal}:
\begin{flalign}
\label{eq:eulerRenorm}
 \E_{T_{(1)},X_{T_{(1)}}}\big[ (\sigma_0(T_{(1)},X^{(1)}_{T_{(1)}})^\top)^{-1} \frac{\hat W^{(1,p),1}_{\min( \delta t,\Delta T_{(1,p)}}}{\min( \delta t,\Delta T_{(1,p)})} \big(\phi\big(T_{(1,p)},X^{(1,p)}_{T_{(1,p)}}\big)  -\phi\big(T_{(1,p)},X^{(1,p^1)}_{T_{(1,p)}}\big) \big)  \big].
\end{flalign}
\begin{Remark}
The renormalization technique described for the renormalization technique of section \ref{subsec:renorm} can be straightforwardly adapted  to the renormalization scheme with antithetics of 
section \ref{subsec:renormAnt}.
\end{Remark}
Of course using equation \eqref{eq:eulerOrig} we expect that variance of the scheme will degrade with the diminution of the time step and we expect the scheme
\eqref{eq:eulerRenorm} to correct this behaviour.
On figure \ref{figEuler} we give the error estimations given by the original scheme and the renormalization technique  (with anithetics of section \ref{subsec:renorm}) depending on the time step for a case with
 burgers non linearity in dimension 4 with $1e6$ particles: as we refine the time step the scheme \eqref{eq:eulerOrig} becomes unusable while the scheme \eqref{eq:eulerRenorm}
gives stable results.

\section{The full non linear case}
In order to  treat some full non linear case, so with a second order derivative $D^2u$ in $f$, the re-normalization technique is necessary as no distribution can meet the finite variance requirement
 even when $f$ is linear in $D^2u$ (see \cite{LOTTW}).\\
Suppose that the $f$ function is as follows :
\begin{flalign*}
f(t,x,y,z,\gamma)
~:=\!
h(t,x)+
\! c(t,x) y^{\ell_0}
\prod_{i=1}^m \big( (b_i \cdot z)^{\ell^1_i} \big)
\prod_{i=m+1}^{2m} \big( (a_i : \gamma)^{\ell_i} \big),
\end{flalign*}
for a given $(\ell_0, \ell_1, \cdot, \ell_m, \ell_{m+1}, \cdots, \ell_{2m}) \in \N^{1+2m}$,
$m \ge 1$, $b_i:[0,T] \x \R^d \to \R^d$ for $i=1, \cdots, m$ are  bounded continuous,
 $h: [0,T] \x \R^d \to \R$  is a  bounded continuous function, and $a_i : [0,T] \x \R^d \to \M^d$, for   $i=m+1, \cdots, 2m $ are  bounded continuous functions.
We note $L= \sum_{i=0}^{2m} \ell_i$.\\
We use a similar algorithm to the one proposed in section \ref{sec:origAlgo}. 
Instead of approximating $f$ using representation \eqref{eq:calU}, we have to take into account  the $D^2u$ term :
\begin{flalign}
\label{eq:calD2U}
[c u^{\ell_0} & \prod_{i=1}^m (b_i \cdot Du)^{\ell^1_i} \prod_{i=m+1}^{2m} (a_i :D^2u)^{\ell_i} ](T_{(1)},X_{T_{(1)}})  = \nonumber\\ &c \prod_{j=1}^{\ell_0} \E_{T_{(1)},X_{T_{(1)}}} \big[ \phi\big(T_{(1,j)},X^{(1)}_{T_{(1,j)}}\big)\big]  \nonumber  \\ & \prod_{i=1}^m ( b_i(T_{(1)},X_{T_{(1)}}).D \E_{T_{(1)},X_{T_{(1)}}}\big[\phi\big(T_{(1,p)},X^{(1,p)}_{T_{(1,p)}}\big)\big])^{\ell_p^1} \nonumber\\ &
 \prod_{i=m+1}^{2m} (a_i :D^2 \E_{T_{(1)},X_{T_{(1)}}}\big[\phi\big(T_{(1,p)},X^{(1,p)}_{T_{(1,p)}}\big)\big])^{\ell_i}.
\end{flalign}
The terms $$\E_{T_{(1)},X_{T_{(1)}}} \big[ \phi\big(T_{(1,j)},X^{(1)}_{T_{(1,j)}}\big)\big]$$ and $$( b_i(T_{(1)},X_{T_{(1)}}).D \E_{T_{(1)},X_{T_{(1)}}}\big[\phi\big(T_{(1,p)},X^{(1,p)}_{T_{(1,p)}}\big)\big])$$ are approximated by the different schemes previously seen.
It remains to give  an approximation of the $(a_i :D^2 \E_{T_{(1)},X_{T_{(1)}}}\big[\phi\big(T_{(1,p)},X^{(1,p)}_{T_{(1,p)}}\big)\big] u)$ term.

\subsection{Ghost particles of dimension $q$}
We extend the definition given in \cite{LTTW} of ghost tree for the full non linear case.
For a particle in dimension $(1)$ of generation $n=1$, we introduce  $q$ associated ghost particles denoted $(1^i)$ for $i=1,...,q$.
Let $ \Kct^1_T := \{(1), (1^{1}), ...,  (1^q) \}$
Then given the collection $\Kct^n_T$ of all particles and ghost particles of generation $n$, 
we define $\Kct^{n+1}_T$ as follows.
Given $k=(k_1, \cdots, k_n) \in \Kct^n_T$, we denote by $o(k)$ its original particle; and when $k_n \in \N$, we denote $k^{i} := (k_1, \cdots, k_{n-1}, k_n^{i})$  for $i \in [1,q]$ and $i$ is noted the order of $k^{i}$. The function $\kappa$ allows us to give the order of a particle  for $k=(k_1, \cdots, k_n) \in \Kct^n_T$ :
\begin{flalign*}
\kappa(k) & =  i,  \mbox{ if } k_n = p^i \mbox{ with } p \in \N, \\
\kappa(k) & =   0, \mbox{ if } k_n = p \mbox{ with } p \in \N,
\end{flalign*}
The variables $T_k$ as well as the mark $\theta_k$ inherits that of the original particle $o(k)$. Similarly $\Delta T_k = \Delta T_{o(k)}$.
Denote also $\Kch^n_T := \{ k \in \Kct^n_T ~: o(k) \in \Kc^n_T \}$.
For every $k = (k_1, \cdots, k_n) \in  \Kct^n_T \setminus \Kch^n_T$, 
we define the collection of its offspring particles by
$$h(k) := \{(k_1, \cdots, k_i, 1), \cdots, (k_1, \cdots, k_i, L) \},$$
and generalizing the definition in section \ref{sec:SemiLinRenorm}, we introduce $q$ collections of all offspring ghost particles:
$$
S^{i}(k)
~~:=~~
\big\{ (k_1, \cdots, k_n, 1^{i}), \cdots, (k_1, \cdots, k_n, L^{i}) \big\},  \mbox{ for } i=1,...,q
$$ 
Then the collection $\Kct^{n+1}_T$ of all particles and ghost particles of generation $n+1$ is given by
$$
\Kct^{n+1}_T
~:=~
\cup_{k \in  \Kct^n_T \setminus \Kch^n_T } \big( S(k) \cup S^{1}(k) \cup ... \cup  S^{q}(k)  \big).
$$
\subsection{$D^2u$ approximations}
In this section, we give some different schemes that can be used  to approximate the $D^2u$ term and that we will compared on some numerical test cases.

\subsubsection{The original $D^2u$ approximation}
The approximation developed in this paragraph was first proposed  in \cite{LTTW} and uses some ghost particle of dimension $q=2$.
To obtained the position of a particle, we freeze its position if its order is $2$ and  inverse its increment if its order is $1$ , so  for every $k = (k_1, \cdots, k_n) \in \Kct^n_T$
\begin{flalign}
\label{eq:brownRenormSecondOrder}
W^{k}_s
~:=~
W^{k-}_{T_{k-}}
		~+~
		\1_{ \kappa(k)=0}
		\hat W^{o(k)}_{s - T_{k-}}  ~-~ \1_{\kappa(k)=1}
		\hat W^{o(k)}_{s - T_{k-}}, \nonumber \\
		~~~\mbox{and}~~
		X^{k}_s := \mu s +\sigma_0  W^k_s,
		~~~\forall s \in [T_{k-}, T_k].
\end{flalign}

Then we use the following representation for the $D^2u$ term in equation \eqref{eq:calD2U} :
\begin{flalign}
\label{eq:secOrderOr}
 D^2  & \E_{T_{(1)},X_{T_{(1)}}}\big[\phi\big(T_{(1,p)},X^{(1,p)}_{T_{(1,p)}}\big)\big]  = \nonumber \\
  & \E_{T_{(1)},X_{T_{(1)}}}\big[ (\sigma_0^{\top})^{-1} \frac{\hat W^{(1,p)}_{\Delta T_{(1,p)}}(\hat W^{(1,p)}_{\Delta T_{(1,p)}})^{\top} - \Delta T_{(1,p)} I_d}{(\Delta T_{(1,p)})^2} \sigma_0^{-1} \psi \big],  
\end{flalign}
where
\begin{flalign*}
\psi= \frac{1}{2} \big[ \phi\big(T_{(1,p)},X^{(1,p)}_{T_{(1,p)}} \big) + \phi\big(T_{(1,p)},X^{(1,p^1)}_{T_{(1,p)}}\big) - 2 \phi\big(T_{(1,p)},X^{(1,p^2)}_{T_{(1,p)}}  \big) \big].
\end{flalign*}
Using for example the equation \eqref{eq:folRep} for the first derivative $Du$, \cite{LTTW} gave the following re-normalized estimator defined by a backward induction:
let $\psih_k := \frac{g(X^k_T)}{\Fb(\Delta T_k)}$ for every $k \in \Kch_T$, then let
\begin{flalign}\label{eq:ghostRepSecOrder}
\psih_k
& ~:=~ 
\frac{1}{\rho(\Delta T_k) } \big( h(T_k, X^k_{T_k})+c(T_k, X^k_{T_k})
\prod_{{\tilde k} \in S(k)} \!\!\! \big(\psih_{{\tilde k}}\1_{\theta({\tilde k})=0} + (\psih_{{\tilde k}} - \psih_{{\tilde k}^2}) \1_{ 1 \le \theta({\tilde k}) \le m } +  \nonumber \\
&  \frac{1}{2}(\psih_{{\tilde k}}+\psih_{{\tilde k}^1} - \psih_{{\tilde k}^2}) \1_{m+1 \le \theta({\tilde k}) \le 2m\}} \big) \Wc_{{\tilde k}} \big),
		~~~~\mbox{for}~k \in \Kct_T \setminus \Kch_T.
\end{flalign}
where 
\begin{flalign*}
\Wc_k 
&:= 
\1_{\{\theta_k = 0\}} 
+ \1_{\{\theta_k  \in \{1, \cdots, m\}\}} \frac{b_{\theta_k}(T_{k-}, X^k_{T_{k-}}) \cdot (\sigma_0^{\top})^{-1} {\hat W}^{o(k)}_{\Delta T_k}}{\Delta T_k}\\
&+~ \1_{\{\theta_k \in \{m+1, \cdots, 2m\}\}} a_{\theta_k} : (\sigma_0^{\top})^{-1} \frac{{\hat W}^{o(k)}_{\Delta T_k} {\hat W}^{o(k)}_{\Delta T_k} - \Delta T_k I_d}{(\Delta T_k)^2} \sigma_0^{-1}.
\end{flalign*}
Then we have 
\begin{flalign*}
u(0,x)
= \E_{0,x}\Big[  \psih_{(1)}\Big].
\end{flalign*}

\subsection{A second representation}

This second representation  uses some ghost particle of dimension $q=3$.
Let $$(\hat W^{k,i})_{k = (k_1, \cdots, k_{n-1}, k_n) \in \N^n, n>1, i=1, 2 }$$ be a sequence of independent $d$-dimensional Brownian motion, which is also independent of $(\Delta T_{k})_{k = (k_1, \cdots, k_{n-1}, k_n) \in \N^n, n>1}$.
The dynamic of the original particles and the ghosts is given by :
\begin{flalign}\label{eq:brownRenormSecondOrderSecRep}
W^{k}_s
~:=~
W^{k-}_{T_{k-}}
		~+~
		\1_{ \kappa(k)=0}
		\frac{\hat W^{o(k),1}_{s - T_{k-}} +\hat W^{o(k),2}_{s - T_{k-}}}{\sqrt{2}}   ~+~ \1_{\kappa(k)=1} \frac{\hat W^{o(k),1}_{s - T_{k-}}}{\sqrt{2}}  ~+~ \1_{\kappa(k)=2} \frac{\hat W^{o(k),2}_{s - T_{k-}}}{\sqrt{2}} \nonumber \\
		~~~\mbox{and}~~
		X^{k}_s := \mu s +\sigma_0  W^k_s,
		~~~\forall s \in [T_{k-}, T_k].
\end{flalign}

We then replace \eqref{eq:secOrderOr} by
\begin{flalign}
\label{eq:secOrderOrSecRep}
 D^2  \E_{T_{(1)},X_{T_{(1)}}}\big[\phi\big(T_{(1,p)},X^{(1,p)}_{T_{(1,p)}}\big)\big]  = \E_{T_{(1)},X_{T_{(1)}}}\big[ 2 (\sigma_0^{\top})^{-1} \frac{\hat W^{(1,p),1}_{\Delta T_k}(\hat W^{(1,p),2}_{\Delta T_k})^{\top}}{(\Delta T_{(1,p)})^2} \sigma_0^{-1} \psi) \big],  
\end{flalign}
where
\begin{flalign*}
\psi = \phi\big(T_{(1,p)},X^{(1,p)}_{T_{(1,p)}} \big) + \phi\big(T_{(1,p)},X^{(1,p^3)}_{T_{(1,p)}}\big) -  \phi\big(T_{(1,p)},X^{(1,p^1)}_{T_{(1,p)}}\big)-  \phi\big(T_{(1,p)},X^{(1,p^2)}_{T_{(1,p)}} \big).
\end{flalign*}
This scheme can be can be easily obtained by applying the differentiation rule used for semi linear equations  on two successive steps with size $\frac{\Delta T_{(1,p)}}{2}$.
A simple calculation  shows that the original scheme has a variance bounded by $|D^2u|_\infty^2 \frac{39}{2}$ while this one has  variance bounded by $|D^2u|_\infty^2 9$ so we expect a diminution of the variance observed with this new  scheme.
\begin{Remark}
This derivation on two consecutive time steps has already been used implicitly for example in \cite{FTW} and already was  numerically superior to a scheme directly using second order Malliavin weight.
\end{Remark}
Recursively the re-normalized estimator is defined by a backward induction:
let $\psih_k := \frac{g(X^k_T)}{\Fb(\Delta T_k)}$ for every $k \in \Kch_T$, then let
\begin{flalign}\label{eq:ghostRepSecOrderSecRep}
\psih_k
& := 
\frac{1}{\rho(\Delta T_k) } \big( h(T_k, X^k_{T_k})+c(T_k, X^k_{T_k})
\prod_{{\tilde k} \in S(k)} \!\!\! \big(\psih_{{\tilde k}}\1_{\theta({\tilde k})=0} + (\psih_{{\tilde k}} - \psih_{{\tilde k}^3}) \1_{ 1 \le \theta({\tilde k}) \le m } + \nonumber \\
&     (\psih_{{\tilde k}}+\psih_{{\tilde k}^3} - \psih_{{\tilde k}^1} - \psih_{{\tilde k}^2}) \1_{m+1 \le \theta({\tilde k}) \le 2m\}} \big) \Wc_{{\tilde k}} \big),
		~~~~\mbox{for}~k \in \Kct_T \setminus \Kch_T.
\end{flalign}
where 
\begin{flalign}
\label{eq:weifhtSecondRep}
\Wc_k 
&:= 
\1_{\{\theta_k = 0\}} 
+ \1_{\{\theta_k  \in \{1, \cdots, m\}\}} \frac{b_{\theta_k}(T_{k-}, X^k_{T_{k-}}) \cdot (\sigma_0^{\top})^{-1} {\hat W}^{o(k),1}_{\Delta T_k}}{\Delta T_k} \nonumber \\
&+~ \1_{\{\theta_k \in \{m+1, \cdots, 2m\}\}} a_{\theta_k} : 2  (\sigma_0^{\top})^{-1} \frac{ {\hat W}^{o(k),1}_{\Delta T_k}  {\hat W}^{o(k),2}_{\Delta T_k}}{(\Delta T_k)^2} \sigma_0^{-1}.
\end{flalign}
Then we have 
\begin{flalign*}
u(0,x)
= \E_{0,x}\Big[  \psih_{(1)}\Big].
\end{flalign*}

\subsection{A third representation}
This representation is only the antithetic version  of the second one and uses some ghost particle of dimension $q=6$.
The dynamic of the original particles and the ghosts is given by :
\begin{flalign}\label{eq:brownRenormSecondOrderThirdRep}
W^{k}_s
& :=
W^{k-}_{T_{k-}}
		~+~
		\1_{ \kappa(k)=0}
		\frac{\hat W^{o(k),1}_{s - T_{k-}} +\hat W^{o(k),2}_{s - T_{k-}}}{\sqrt{2}}   ~+~ \1_{\kappa(k)=1} \frac{\hat W^{o(k),1}_{s - T_{k-}}}{\sqrt{2}}  ~+~ \1_{\kappa(k)=2} \frac{\hat W^{o(k),2}_{s - T_{k-}}}{\sqrt{2}} ~-~\nonumber \\
&  \1_{\kappa(k)=4}\frac{\hat W^{o(k),1}_{s - T_{k-}} +\hat W^{o(k),2}_{s - T_{k-}}}{\sqrt{2}}  ~-~\1_{\kappa(k)=5} \frac{\hat W^{o(k),1}_{s - T_{k-}}}{\sqrt{2}} ~-~ \1_{\kappa(k)=6} \frac{\hat W^{o(k),2}_{s - T_{k-}}}{\sqrt{2}} \nonumber \\
&		~~~\mbox{and}~~
		X^{k}_s := \mu s +\sigma_0  W^k_s,
		~~~\forall s \in [T_{k-}, T_k].
\end{flalign}

We then replace \eqref{eq:secOrderOr} by
\begin{flalign}
\label{eq:secOrderOrSecRep}
 D^2 &  \E_{T_{(1)},X_{T_{(1)}}}\big[\phi\big(T_{(1,p)},X^{(1,p)}_{T_{(1,p)}}\big)\big]  = \nonumber \\
    &  \E_{T_{(1)},X_{T_{(1)}}}\big[  (\sigma_0^{\top})^{-1} \frac{\hat W^{(1,p),1}(\hat W^{(1,p),2})^{\top}}{(\Delta T_{(1,p)})^2} \sigma_0^{-1} \psi) \big],  
\end{flalign}
where
\begin{flalign*}
\psi &  =   \phi\big(T_{(1,p)},X^{(1,p)}_{T_{(1,p)}} \big) + 2 \phi\big(T_{(1,p)},X^{(1,p^3)}_{T_{(1,p)}}\big) -  \phi\big(T_{(1,p)},X^{(1,p^1)}_{T_{(1,p)}}\big)-  \phi\big(T_{(1,p)},X^{(1,p^2)}_{T_{(1,p)}} \big) + \\
&   \phi\big(T_{(1,p)},X^{(1,p)}_{T_{(1,p^4)}} \big)- \phi\big(T_{(1,p)},X^{(1,p^5)}_{T_{(1,p)}}\big)- \phi\big(T_{(1,p)},X^{(1,p^6)}_{T_{(1,p)}}\big) \big),
\end{flalign*}
and the weights are still given by equation \eqref{eq:weifhtSecondRep}.
The backward induction is defined as follows:
let $\psih_k := \frac{g(X^k_T)}{\Fb(\Delta T_k)}$ for every $k \in \Kch_T$, then let
\begin{flalign} \label{eq:ghostRepSecOrderSecRep3}
\psih_k
& :=
\frac{1}{\rho(\Delta T_k) } \big( h(T_k, X^k_{T_k})+ \frac{c(T_k, X^k_{T_k})}{2}
\prod_{{\tilde k} \in S(k)} \!\!\! \big((\psih_{{\tilde k}}+\psih_{{\tilde k}^4})\1_{\theta({\tilde k})=0} + (\psih_{{\tilde k}} - \psih_{{\tilde k}^4}) \1_{ 1 \le \theta({\tilde k}) \le m } + \nonumber \\
&     \frac{1}{2} (\psih_{{\tilde k}}+ 2\psih_{{\tilde k}^3} - \psih_{{\tilde k}^1} - \psih_{{\tilde k}^2} +\psih_{{\tilde k}^4}-  \psih_{{\tilde k}^5} - \psih_{{\tilde k}^6}) \1_{m+1 \le \theta({\tilde k}) \le 2m\}} \big) \Wc_{{\tilde k}} \big),
		~~~~\mbox{for}~k \in \Kct_T \setminus \Kch_T.
\end{flalign}
where the weights are given by equation \eqref{eq:weifhtSecondRep}.
And as usual we have 
\begin{flalign*}
u(0,x)
= \E_{0,x}\Big[  \psih_{(1)}\Big].
\end{flalign*}

\begin{Remark}
Extension to schemes for derivatives of order more than 3 is obvious with the two last schemes. 
\end{Remark}
\subsection{Numerical results}
For all test cases in this section we take $\mu = 0.2 \1 $, $\sigma_0 = 0.5 \un$ and we want to evaluate $u(0,0.5 \1)$.
We test the 3 schemes previously described :
\begin{itemize}
\item Version 1 stands for the original version of the scheme using backward recursion \eqref{eq:ghostRepSecOrder},
\item Version 2 stands for the second representation using second backward recursion \eqref{eq:ghostRepSecOrderSecRep},
\item Version 3 stands  for the third representation corresponding to the antithetic version of the second representation and using backward recursion 
\eqref{eq:ghostRepSecOrderSecRep3}. Notice that in this case all terms in $u$ in $f$ are treated with antithetic ghosts.
\end{itemize} 
We give results for the non nested version as the nested version doesn't improve the results very much.
\begin{itemize}
\item  We first choose a non linearity  $$f(u,Du,D^2u) = h(t,x)+\frac{0.1}{d} u(\un:D^2u),$$ where  $\mu = 0.2 \1 $, $\sigma_0 = 0.5 \un$ and
\begin{flalign*}
h(t,x)= &(\alpha+ \frac{\sigma_0^2}{2}) \cos(x_1+..+x_d) e^{\alpha(T-t)}+ \\
&  0.1 \cos(x_1+..+x_d)^2 e^{2\alpha(T-t)}+ \mu \sin(x_1+..+x_d) e^{\alpha(T-t)},
\end{flalign*}
 with $\alpha=0.2$.
We suppose that the final solution is given by $g(x)= \cos(x_1+..+x_d)$ such that the analytical solution is 
\begin{flalign*}
u(t,x)= \cos(x_1+..+x_d) e^{\alpha(T-t)}.
\end{flalign*}
This test case will be noted test C.
In the example we want evaluate $u(0,0.5 \1)$.
 First we take $d=4$ and give   the results obtained for different maturities on figures \ref{figNL} and \ref{figNL1}.
\begin{figure}[H]
  \centering
  \includegraphics[width=7cm]{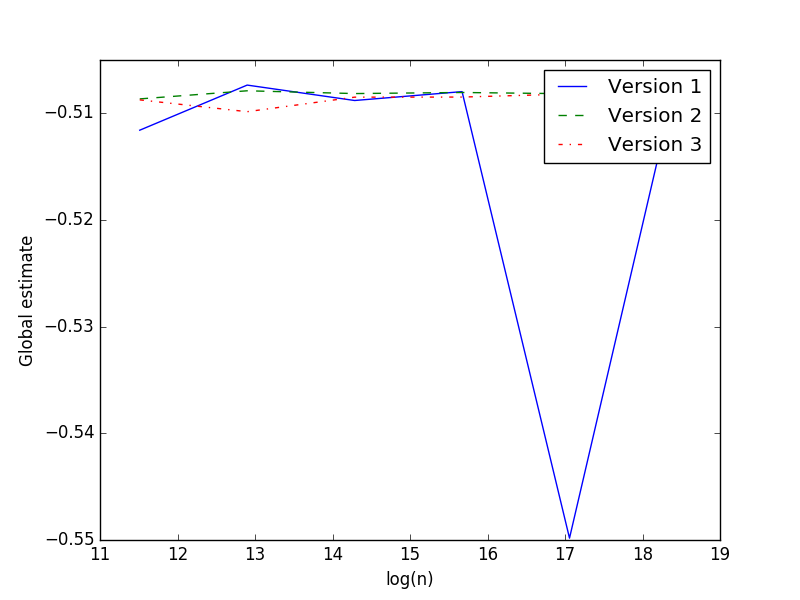}
  \includegraphics[width=7cm]{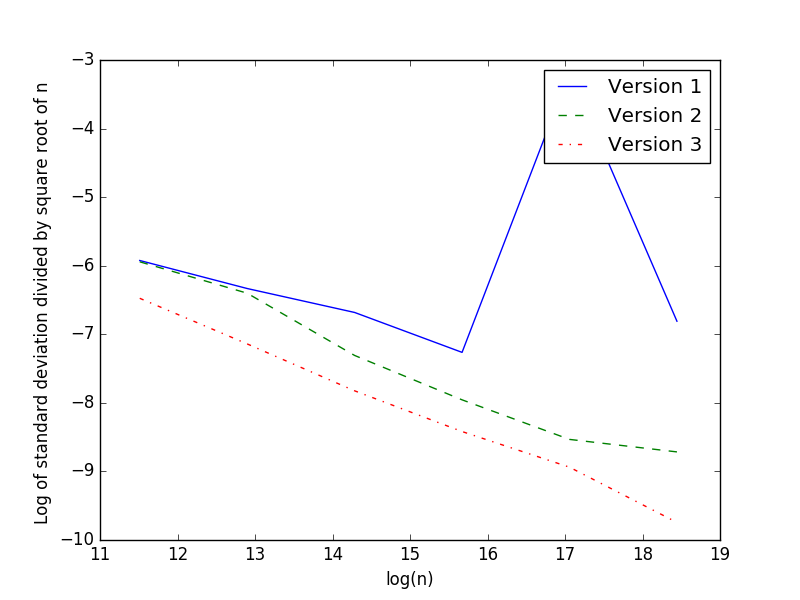}
   \caption{Solution and error obtained  in $d=4$  for test case C with $T=1$, analytical solution is $-0.50828$.}
\label{figNL}
\end{figure}
\begin{figure}[!htb]
  \centering
  \includegraphics[width=7cm]{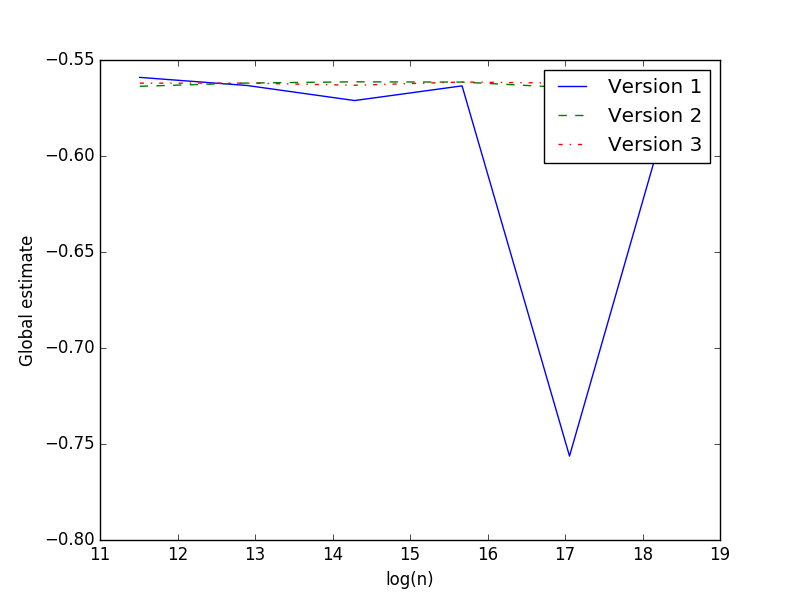}
  \includegraphics[width=7cm]{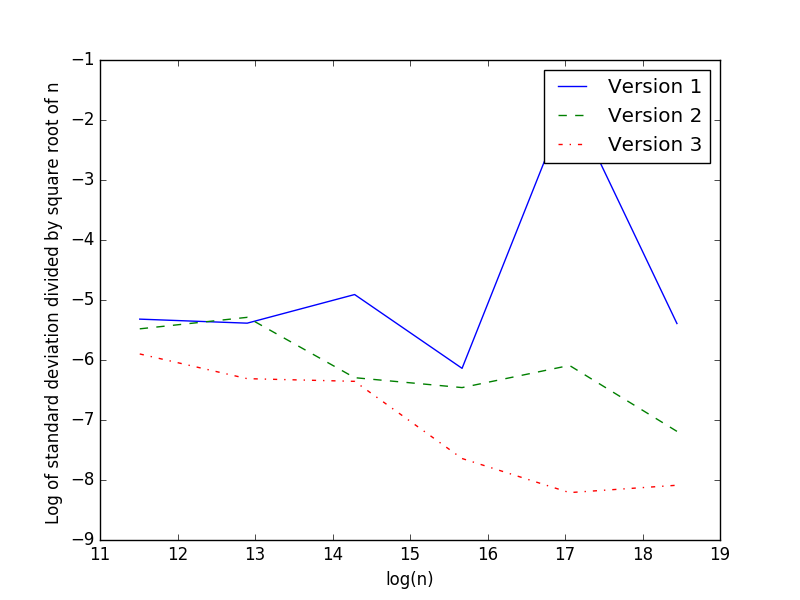}
   \caption{Solution  and error obtained in $d=4$  for test case C with $T=1.5$, analytical solution is $-0.561739$.}
\label{figNL1}
\end{figure}
We then test in dimension 6 the different schemes on figure \ref{figNL2}. Besides on figure \ref{figNL2_} we show that
the schemes provide a good accuracy for the computation of the derivatives by plotting $(\1.Du)$ for the three versions : as expected the
accuracy is however slightly less good than for  the function evaluation.
\begin{figure}[!htb]
  \centering
  \includegraphics[width=7cm]{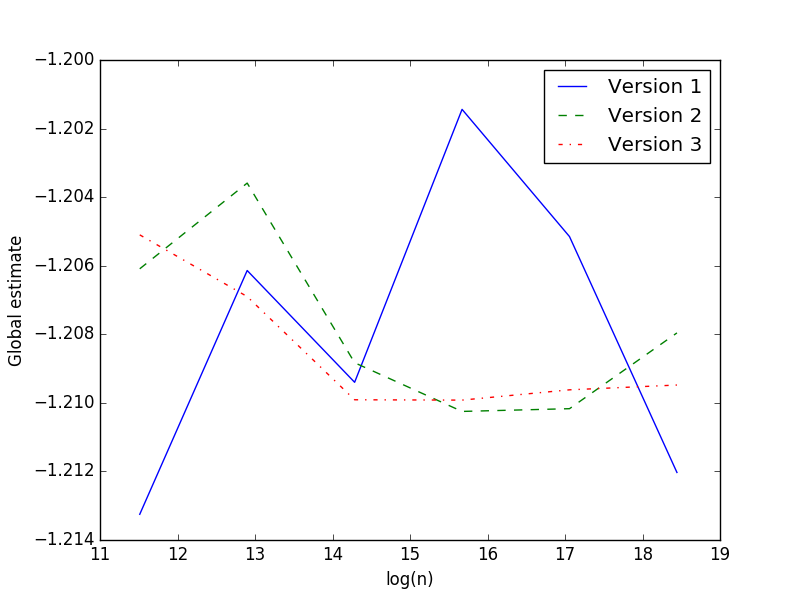}
  \includegraphics[width=7cm]{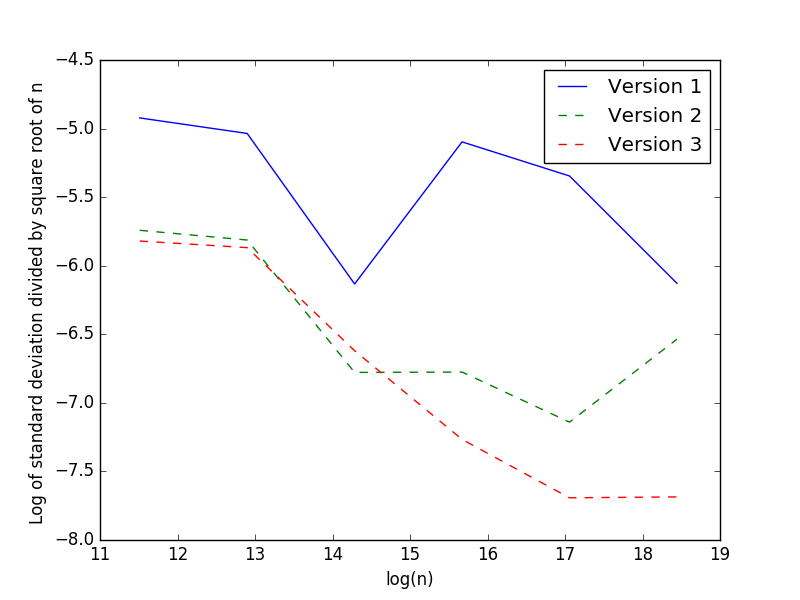}
   \caption{Solution obtained and error in $d=6$ for test case C  with $T=1$, analytical solution is $-1.20918$.}
\label{figNL2}
\end{figure}
\begin{figure}[!htb]
  \centering
  \includegraphics[width=7cm]{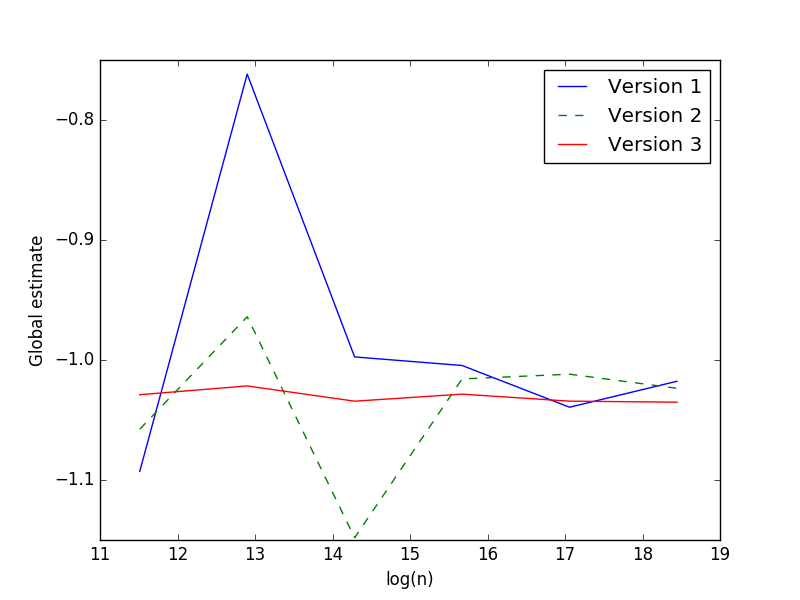}
  \includegraphics[width=7cm]{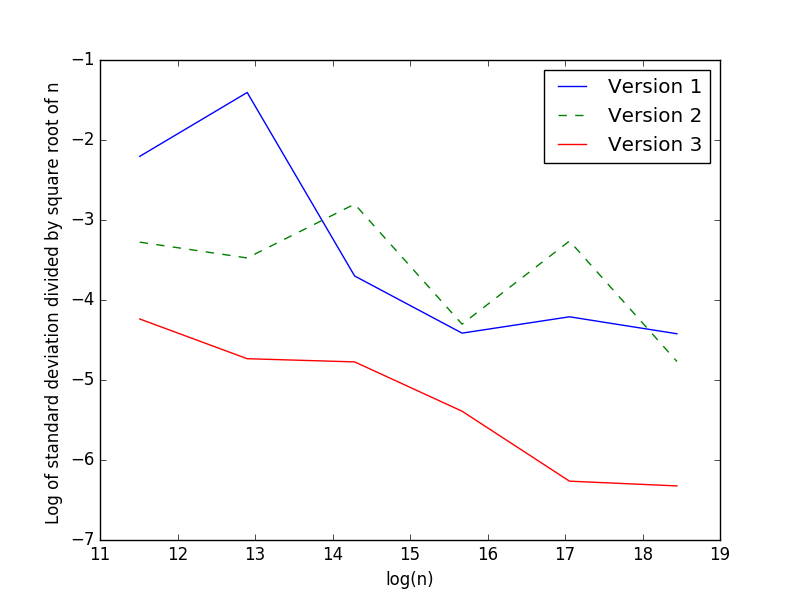}
   \caption{Derivative $(\1.Du)$ obtained and error in $d=6$ for test case C  with $T=1$.}
\label{figNL2_}
\end{figure}
\item At last we consider the test D where $d=4$,  and 
\begin{flalign*}
f(u,Du,D^2u) = 0.0125 (\1.DU) (\un:D^2u).
\end{flalign*}
We give the solution and error obtained for the 3 methods on figure \ref{figNL3}.
\begin{figure}[!htb]
  \centering
  \includegraphics[width=7cm]{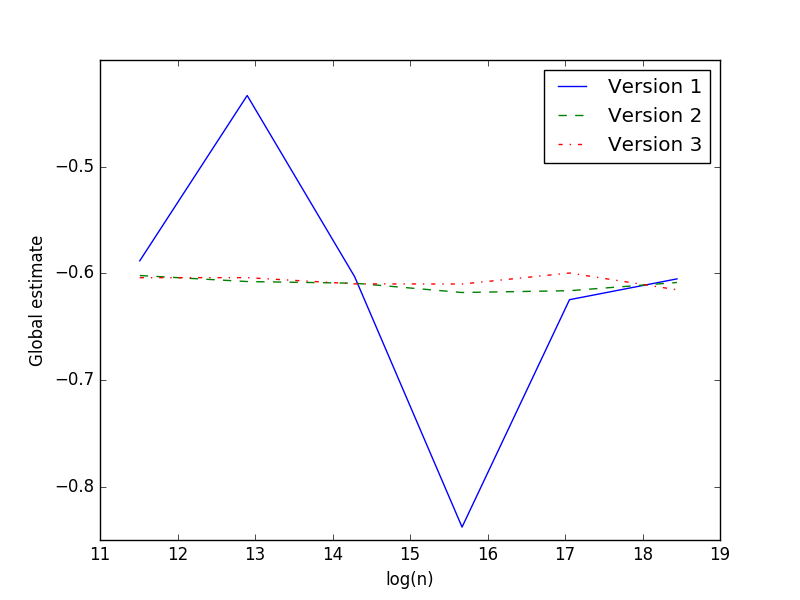}
  \includegraphics[width=7cm]{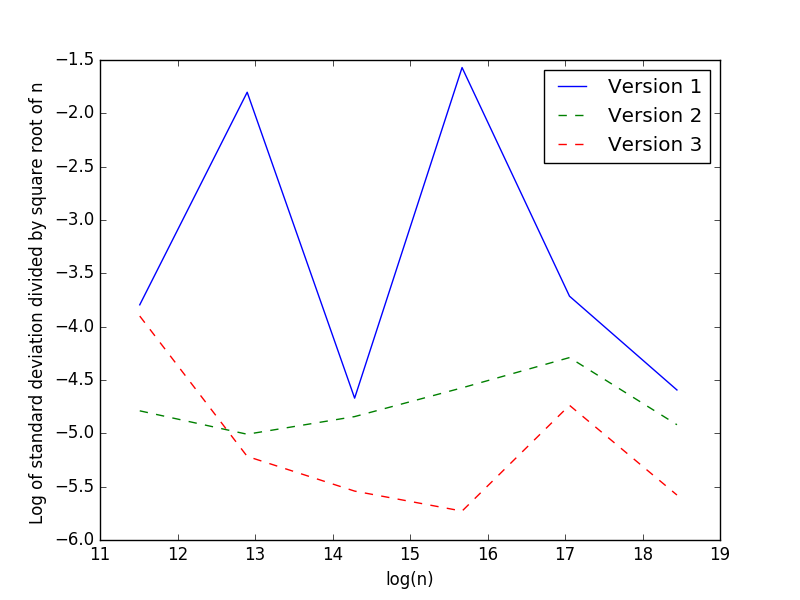}
   \caption{Solution and error obtained for  $d=4$ for test case D with  $T=1$.}
\label{figNL3}
\end{figure}
\end{itemize}
On all the test cases, the last representation using antithetic variables gives the best result in term of variance reduction but at a price of memory consumption increase: as order of the ghost representation increase so does the memory needed.

\section{Conclusion}
As the scheme and methods developped here let us extend the maturities than can be used to evaluate the solution of some semi linear and full non linear equation.
This is achieved by an increase of the computational time and the memory consumption.

\bibliographystyle{plain}
\bibliography{bibVariationBranching}

\begin{thebibliography}{10}

\bibitem{BET}
Bruno Bouchard, Ivar Ekeland, and Nizar Touzi.
\newblock On the malliavin approach to monte carlo approximation of conditional
  expectations.
\newblock {\em Finance and Stochastics}, 8(1):45--71, 2004.

\bibitem{BouchardTouzi}
Bruno Bouchard and Nizar Touzi.
\newblock Discrete-time approximation and monte-carlo simulation of backward
  stochastic differential equations.
\newblock {\em Stochastic Processes and their applications}, 111(2):175--206,
  2004.

\bibitem{BW}
Bruno Bouchard and Xavier Warin.
\newblock Monte-carlo valuation of american options: facts and new algorithms
  to improve existing methods.
\newblock In {\em Numerical methods in finance}, pages 215--255. Springer,
  2012.

\bibitem{cheridito}
Patrick Cheridito, H~Mete Soner, Nizar Touzi, and Nicolas Victoir.
\newblock Second-order backward stochastic differential equations and fully
  nonlinear parabolic pdes.
\newblock {\em Communications on Pure and Applied Mathematics},
  60(7):1081--1110, 2007.

\bibitem{DOW}
Mahamadou Doumbia, Nadia Oudjane, and Xavier Warin.
\newblock Unbiased monte carlo estimate of stochastic differential equations
  expectations.
\newblock {\em to appear in ESAIM P\&S}, 2017.

\bibitem{FTW}
Arash Fahim, Nizar Touzi, and Xavier Warin.
\newblock A probabilistic numerical method for fully nonlinear parabolic pdes.
\newblock {\em The Annals of Applied Probability}, pages 1322--1364, 2011.

\bibitem{FLLLT}
Eric Fourni{\'e}, Jean-Michel Lasry, J{\'e}r{\^o}me Lebuchoux, Pierre-Louis
  Lions, and Nizar Touzi.
\newblock Applications of malliavin calculus to monte carlo methods in finance.
\newblock {\em Finance and Stochastics}, 3(4):391--412, 1999.

\bibitem{LGW1}
Emmanuel Gobet, Jean-Philippe Lemor, and Xavier Warin.
\newblock A regression-based monte carlo method to solve backward stochastic
  differential equations.
\newblock {\em The Annals of Applied Probability}, 15(3):2172--2202, 2005.

\bibitem{LOTTW}
Pierre Henry-Labordere, Nadia Oudjane, Xiaolu Tan, Nizar Touzi, and Xavier
  Warin.
\newblock Branching diffusion representation of semilinear pdes and monte carlo
  approximation.
\newblock {\em arXiv preprint arXiv:1603.01727}, 2016.

\bibitem{HTT2}
Pierre Henry-Labordere, Xiaolu Tan, and Nizar Touzi.
\newblock Unbiased simulation of stochastic differential equations.
\newblock {\em arXiv preprint arXiv:1504.06107}, 2015.

\bibitem{LTTW}
Pierre~Henri Labord\`ere, Xiaolu Tan, Nizar Touzi, and Xavier Warin.
\newblock Truncation and renormalization techniques for solving the nonlinear
  pdes by branching processes, 2017.

\bibitem{LGW2}
Jean-Philippe Lemor, Emmanuel Gobet, and Xavier Warin.
\newblock Rate of convergence of an empirical regression method for solving
  generalized backward stochastic differential equations.
\newblock {\em Bernoulli}, 12(5):889--916, 2006.

\bibitem{PardouxPeng}
Etienne Pardoux and Shige Peng.
\newblock Adapted solution of a backward stochastic differential equation.
\newblock {\em Systems \& Control Letters}, 14(1):55--61, 1990.

\bibitem{zhang}
Jianfeng Zhang.
\newblock A numerical scheme for bsdes.
\newblock {\em the annals of applied probability}, 14(1):459--488, 2004.

\end{thebibliography}

\end{document}